\documentclass[11pt]{article}

\usepackage{amssymb}
\usepackage{latexsym}
\usepackage[mathscr]{euscript}
\usepackage{graphicx}

\usepackage{amscd}
\usepackage{amsmath}
\usepackage{amssymb}
\usepackage{fullpage}

\usepackage{amssymb}
\usepackage{latexsym}
\usepackage[mathscr]{euscript}
\usepackage{pb-diagram}
\usepackage{latexsym}
\usepackage{amsfonts}
\usepackage{epsfig}

\newcommand\commentout[1]{\marginpar{\tiny $\backslash$commentout}}

\newcommand\qed{\hfill$\square$}
\def\column#1#2{\mathrel{\mathop{#1}\limits_{#2}}}

\def\compcirc {\mbox{\hspace{.05cm}}\raisebox{.04cm}{\tiny  {$\circ$ }}}

\newtheorem{Lemma}{Lemma}[section]
\newtheorem{Theorem}{Theorem}[section]
\newtheorem{Proposition}{Proposition}[section]
\newtheorem{Definition}{Definition}[section]
\newtheorem{Corollary}{Corollary}[section]

\newenvironment{Proof}{\par\noindent\textbf{Proof:}}
{\qed}

\usepackage{graphics}
\title{Project Description}

\title{$K$-theoretic descent and a motivic Atiyah-Segal theorem}
\author{Gunnar Carlsson \footnote{Research supported in part by NSF DMS-0406992}  \\ Department of Mathematics, Stanford University  \\Stanford, California 94305}
\begin{document}
\maketitle
\section{Introduction} 
In recent years there have been striking developments in the study of the algebraic $K$-theory of fields (see \cite{voevodsky}, \cite{friedlander}, \cite{survey}).  These developments have resulted in the identification of the $E_2$-term of a spectral sequence converging to the algebraic $K$-theory of fields (rather a mod-$p$ version of it) in terms of Galois cohomology, and have resolved long standing conjectures of Milnor and Bloch-Kato in the area.    In particular the $E_2$-term depends only on the separable Galois  group of the field.  These results raise a pair of questions.   
\begin{itemize}
\item{For a field $F$, is there a construction depending only on the separable Galois  group $G_F$  of $F$, which reconstructs the full $K$-theory spectrum of $F$, rather than an $E_2$-term of a spectral sequence?  Can such a construction be made functorial for inclusions of subfields and the corresponding inclusions of separable Galois  groups? }
\item{The conjectures of Milnor and Bloch-Kato are cohomological statements about the separable Galois  group, specifically that the cohomology is generated in degree one and  that the relations are quadratic.  Are there properties more closely connected to the structure of the separable Galois  group as a group which are relevant?  }
\end{itemize} 

A conjectural answer to the first question in the case of geometric fields (i.e. fields containing an separably closed subfield $k$)  was proposed in \cite{repassembly}, and proved in the case of abelian separable Galois  groups.   The proofs  in \cite{repassembly}  relied on an a priori  explicit understanding of the algebraic $K$-theory of $F$, as well as the complex representation theory of $G_F$, and therefore could not directly lead to a proof of the conjectures in \cite{repassembly}.   The contributions of the present paper are now as follows.  

\begin{enumerate}
\item{The observation that absolute Galois groups of geometric fields share the property of {\em total torsion freeness}, together with an analysis of some simple properties implied by this property.  This property is far from generic among torsion free groups, as examples readily show.  The property has interesting implications for the representation theory of the group.}
\item{Proof that the conjectures in \cite{repassembly} can be interpreted as the statement that a map from a fixed  point spectrum to a mapping spectrum from a certain equivariant pro-scheme into the $K$-theory spectrum is an equivalence.    So, an analogue of  homotopy fixed points are computed in the $\Bbb{A}^1$ category rather than in the usual homotopy category of spaces. This computation uses in an essential way the totally torsion freeness of the absolute Galois groups.  }
\item{As a corollary to point 2 above, one finds that the conjectures in \cite{repassembly} can now be viewed as statements about the invariance of homotopy fixed point sets in the $\Bbb{A}^1$ category.  This leads to a strategy for proving these conjectures which we expect to explore in future work.  } 
\end{enumerate}

We now elaborate a bit on this discussion.  
To describe the program in broad terms, we recall from \cite{repassembly} that there is a map of spectra 
$$ \alpha _F : K(Rep_k[G_F])^{\wedge}_{\epsilon} \rightarrow K(F)^{\wedge}_p
$$
called the {\em representational assembly} 
where $G_F$ denotes the separable Galois  group of a geometric field $F$, $Rep[G_F]$ denotes the category of finite dimensional continuous representations of the profinite group $G_F$ over the algebraically closed ground field $k \subseteq F$ , and where $(-)^{\wedge}_{\epsilon}$ denotes a certain derived completion construction, the technical aspects of which were studied  in \cite{completion}.  The conjecture is that this map is an equivalence for geometric fields.  A difficulty in approaching this problem is that the $K(Rep[G_F])^{\wedge}_{\epsilon}$  does not  share the same good formal properties (e.g. localization and devissage) which the $K$-theory construction enjoys, which makes it difficult to begin any kind of inductions.    In this paper, we restrict ourselves to the case where $G_F$ is  a  pro-$p$ group, where $p$ is prime to the characteristic of $F$, and  identify $K(Rep[G_F])^{\wedge}_{\epsilon}$  with $K({\cal B}G_F)^{\wedge}_p$,  where ${\cal B}G_F$ is a pro-scheme which plays a role in $\Bbb{A}^1$-homotopy theory similar to the role played in ordinary homotopy theory by the classifying space  $BG_F$ (which, since $G_F$ is a profinite group, is a pro-simplicial set).   More generally, if $W$ is a smooth scheme over $Spec(k)$, where $k$ is a an algebraically closed field,  with action by $G_F$, we consider the equivariant $K$-theory $K^{G_F}(W)$, which is a module over the commutative ring spectrum  $K(Rep_k[G_F])$, and show that (Theorem \ref{mainone} and Corollary \ref{maintwo})

\begin{equation} \label{main} K^{G_F}(W) ^{\wedge}_{\epsilon} \cong K({\cal E}G_F \column{\times}{G_F} W)^{\wedge}_p
 \cong K^{G_F} ({\cal E}G_F \times W)^{\wedge}_p
\end{equation} 
where ${\cal E}G_F$ is a $G_F$-equivariant pro-scheme whose orbit pro-scheme is ${\cal B}G_F$.  The $K$-theory  functor is extended to pro-schemes by setting it equal to the corresponding homotopy colimit of spectra.  Note that in the case $W = Spec(k)$, $K^{G_F}(W) $ is the category of finite dimensional representations of $G_F$, and so this simply the map $\alpha _F$ mentioned above.  
We note that this equivalence can only be proved for totally torsion free pro-$p$ groups, and also that due to the representability of algebraic $K$-theory in ${\Bbb{A}}^1$-homotopy theory, the right hand side can be viewed as a $G_F$- homotopy fixed point spectrum of $K(W)$  computed in the $\Bbb{A}^1$ category.  
 If we consider the case $W = Spec(\overline{F})$, where $\overline{F}$ denotes the algebraic closure of $F$, it can be shown that 

\begin{equation} \label{descentequivalence} 
K^{G_F}(Spec(\overline{F}) )\cong K(F)\end{equation} 

 and further  that 
$K(F)^{\wedge}_{\epsilon}  \cong K(F)^{\wedge}_p$, where $K(F)$ is given a $K(Rep_k[G_F])$-module structure via the equivalence (\ref{descentequivalence})  above.  It follows that the equivariant map $Spec(\overline{F}) \rightarrow Spec(k)$ induces a map 
$$  \beta _F: K({\cal B}G_F)^{\wedge}_p \rightarrow K(F)^{\wedge}_p
$$
which after making the identifications described above is identified with the representational assembly $\alpha _F$.  Therefore, in order to prove that $\alpha _F$ is an equivalence, it will suffice to prove that $\beta _F$ is an equivalence.  Since $\beta _F$ can also be identified with the map 

$$  K^{G_F}({\cal E}G_F )  =  K^{G_F}({\cal E}G_F \times Spec(k) ) \rightarrow K^{G_F}({\cal E}G_F \times Spec(\overline{F})) 
$$

it is clear that what is now required is a rigidity statement for the functor $K^{G_F}({\cal E}G_F \times - )^{\wedge}_p$ or equivalently the functor $K({\cal E}G_F \column{\times}{G_F} - )^{\wedge}_p$.    We will return to this rigidity question in a future paper.  

Because of the representability of algebraic $K$-theory in the $\Bbb{A}^1$ homotopy theoretic setting, $K({\cal B}G)$ can be interpreted as a mapping spectrum  $F({\cal B}G, K(k))$.  This suggests an analogy between our results and those of Atiyah \cite{atiyah} and Atiyah-Segal \cite{atiyahsegal} on the mapping spectrum of the classifying space of a finite group, or even a compact Lie group into the periodic $K$-theory spectrum.  Their results are stated in terms of completions of homotopy groups as modules over the representation rings of the finite or compact Lie groups in question, but it is not difficult to see that they can be reinterpreted as space level statements about derived completions over the equivariant $K$-theory spectrum.  Our theorem, then, can be stated as the precise analogue for the Atiyah-Segal theorem with periodic $K$-theory replaced by algebraic $K$-theory and compact Lie groups replaced by profinite groups which are totally torsion free.  It is interesting to speculate about how one might generalize the results to profinite groups with the totally torsion free hypothesis removed.  One possibility is that the space ${\cal E}G$ might no longer be free, but might instead resemble the universal space for the family of finite subgroups in the $G$.  We hope to return to these questions as well in a future paper.  

The author wishes to thank P. Diaconis, D. Dugger, L. Hesselholt, T. Lawson, G. Lyo,  I. Madsen, D. Ramras, R. Vakil, and K. Wickelgren for numerous helpful conversations on this work.  

\section{Preliminaries }
\subsection{ General}

We will assume that the reader is familiar with the theory of ring spectra and modules over them, as well as the notion of a commutative ring spectrum. We will use the version of this theory which uses symmetric spectra, which is presented in \cite{schwede}.  In particular, we will assume the reader is familiar with the standard choices of model structures on this category, and the associated notions of cofibrant objects, as well as $Hom_A(M,N)$, $M \column{\wedge}{A} N$ and universal coefficient and K\"{u}nneth spectral sequences converging to them.    We will also assume familiarity with homotopy limits and colimits in this category.  Throughout the paper, a spectrum will mean a symmetric spectrum, a ring spectrum will mean a symmetric ring spectrum, and a commutative ring spectrum will mean a commutative symmetric ring spectrum.  

We will also require the notion of a {\em special pro-object} and a {\em special $Ind$-object}.  We let $\Bbb{N}$ denote the set of non-negative integers, regarded as an ordered set using the linear ordering, and therefore also as a category in which there is a unique morphism from $m$ to $n$ if and only if $m \geq n$, and so that $Mor(m,n) = \emptyset$ if $m < n$. By a {\em special pro-object} in a category $\underline{C}$ we will mean a functor $\Bbb{N} \rightarrow \underline{C}$, and by a morphism of special pro-objects we will mean a natural transformation of functors.  The category of special pro-objects in $\underline{C}$ will be denoted by $\mbox{sPro}\underline{C}$.  A special $Ind$-object and a morphism of such are defined in the dual way.

\subsection{$K$-theory} \label{ktheory}

In this paper, we will use the Waldhausen version of $K$-theory, constructed in \cite{waldhausen}  using his $S_.$-construction.  This construction takes as its input a  {\em category with cofibrations and weak equivalences}  or a {\em Waldhausen category}, which is more general than the input required for Quillen's original $BQ$ construction in \cite{quillen}.  Waldhausen shows that the exact categories of Quillen can be used to construct a Waldhausen category, and that when the $S_.$ construction is applied to this Waldhausen category, the result is equivalent to Quillen's construction.  The advantages of the Waldhausen construction are that it simultaneously carries the structure of a symmetric spectrum and  a  ring spectrum,  and moreover admits devissage and localization arguments.  We state  the results we will need in this paper.  

\begin{Proposition} \label{commutative} Let $C$ denote a symmetric  monoidal abelian category, which is strict as a monoidal category.  Let $K(C)$ denote the $S_.$-construction on the corresponding Waldhausen category.  Then $K(C)$ is a naturally a commutative  ring spectrum. 
\end{Proposition} 
\begin{Proof}  The fact that $K(C)$ is a  ring spectrum is proved in \cite{geisser}.  The commutativity result  is a direct extension of  the methods of \cite{geisser}.  
\end{Proof} 

\begin{Proposition}  \label{modulestructure} Let $C \hookrightarrow D$ be an inclusion of symmetric monoidal abelian categories, both of which are strict as monoidal categories.  Suppose further that we are given an abelian subcategory intermediate between $C$ and $D$, so we have $C \subseteq E \subseteq D$, and that the monoidal pairing of $D$, when restricted to $C \times E$ has its image in $E$.  Then the  spectrum $K(E)$ is naturally a module over the commutative  ring spectrum $K(C)$, and further any inclusion $E \rightarrow E^{\prime}$ of such subcategories induces  a homomorphism of $K(C)$-module spectra.  
\end{Proposition} 
\begin{Proof} Follows directly from the methods of \cite{geisser}.  
\end{Proof}

\begin{Proposition}  Let $A \subseteq B$ be a Serre subcategory in an abelian category $B$, and let the quotient category $B/A$ be defined as in \cite{swan}.  Then there is a homotopy fibration sequence of spectra 

$$  K(A) \rightarrow K(B) \rightarrow K(B/A)
$$
Moreover, if $A $ is symmetric monoidal, and strict as a monoidal abelian category, then $K(B)$ and $K(B/A) $ are  in a natural way  module spectra over the commutative  ring spectrum $K(A)$, and the map $K(B) \rightarrow K(B/A) $ induced by the inclusion $B \hookrightarrow B/A$ is a morphism of $K(A)$-module spectra.   
\end{Proposition}
\begin{Proof} This is Waldhausen's localization theorem in \cite{waldhausen}, together with a straightforward application of the methods of \cite{geisser}.  
\end{Proof}

\subsection{Equivariant algebraic $K$-theory} \label{equivariant}
We will fix a separably closed ground field $k$ throughout this paper.  We will also occasionally  use the notation $*$ to denote $Spec(k)$.   Let $X$ denote a Noetherian scheme over $Spec(k)$,   on which a finite group $G$ acts by automorphisms.  The action is given by a homomorphism $\rho: G \rightarrow Aut (X)$. 

\begin{Definition} By a {\em $G$-twisted $X$-module}, we will mean a coherent  $X$-module $M$equipped with a $k$-linear $G$-action so that $g\cdot (\alpha m) = \alpha ^g (g\cdot m)$ for all $g  \in G, \alpha \in  {\cal O}_X,  $ and $m \in M$.  $G$-twisted $X$-modules form a category whose morphisms are the equivariant $X$-linear maps, and in fact admits the structure of an exact category in the sense of Quillen \cite{quillen} by restriction of the usual structure from the category of all $A$-modules without $G$-action.  We denote this category by $\mbox{\em Mod}^G(X)$.  

\end{Definition} 


One can now consider the algebraic $K$-theory of this category, and denote it by $K^G(X)$.  Note that this means that we are considering what is usually referred to as $G$-theory  or $K^{\prime}$-theory.  Since all schemes we will be considering will be smooth, this distinction will not affect the $K$-theory spectra, but will make computational devices such as localization more directly available.  In the case of a point, i.e. $X = Spec(k)$, we have that $K^G(*)$ is simply the category of finite dimensional $k$-linear representations of $G$.  $\pi _0 (K^G(k))$ is then the $k$-linear representation ring $R_k[G]$.  In particular, if  either (1) $Char(k) = 0$ or (2) $Char(k)$ is relatively prime to $\# (G)$, then $\pi _0K^G(*)$ is isomorphic to the usual complex representation ring or character ring  of $G$.  $K^G(*)$ is a commutative  ring spectrum.     It is also equipped with an augmentation $\epsilon : K^G(*) \rightarrow K(*)$, which forgets the $G$ action, and for which applying $\pi _0$ induces the usual augmentation on representation rings, which we also refer to as $\epsilon$. Let $I_k[G]$ denote the augmentation ideal, i.e the kernel of $\epsilon$.  
 We consider any Noetherian scheme $X$ over $Spec(k)$  with $G$ action, and let $\pi _0 (X)$ denote the (finite) set of connected components of $X$.  The group $G$ acts on the set $\pi _0 (X)$, and we let $\{ [x_1], \ldots , [x_n] \}$ denote an orbit decomposition of $\pi _0 (X)$, with each $x_i$ an orbit representative. Let $G_i$ denote the stabilizer of the component $x_i$,  let $X_i$ denote the subscheme consisting of points belonging to the components in the orbit $[x_i]$, and let $\overline{X}_i$ denote the set of points belonging to the component $x_i$.  $X_i$ is a $G$-subscheme, and $\overline{X}_i$ is a $G_i$-subscheme.  The we have the following result. 
 
 \begin{Proposition} \label{disjoint} 
 There are equivalences of $K^G(*)$-module spectra
 $$ K^G(X) \simeq \prod _{i=1}^n K^G(X_i)  \simeq \prod_{i=1}^n K^{G_i}(\overline{X}_i)
 $$
 \end{Proposition} 
 
 We recall that for the  action of a finite group  $G$ on a scheme $X$, there is always an orbit scheme $X/G$ (see \cite{mumford}).  We have the following result in the special case of a free $G$-action. 
 
 \begin{Proposition}  \label{free} Suppose that the action of $G$ on $X$ is free in the sense that it is free on the set  of all points,  where a point denotes a morphism $Spec(F) \rightarrow X$  for some field $F$ containing $k$.  Then there is a natural equivalence $K^G(X) \cong K(X/G)$.  
 \end{Proposition} 
 \begin{Proof}  The equivalence of categories from which this equivalence results is proved in \cite{grothendieck}.   \end{Proof} 
 
 We also have the following more specialized result. 
 
 \begin{Proposition}  Let $X$ be any  free $G$-$k$-scheme.  Suppose that $X$ is smooth, of finite dimension $d$.  Then the ideal $I[G]^{d+1}_k$ annihilates the $\pi _0 K^G(*)$-modules $\pi _i K^G(X)$.  
 \end{Proposition} 
 \begin{Proof} We consider first the case of a field extension over $k$.  In this case, we are given an action of the group $G$ on $F$, with fixed point subfield $F^G$, and the $K$-theory of the associated category of descent data is equivalent to the category of finite dimensional $F^G$ vector spaces.  In particular, $\pi _0 K^G(Spec(F))$ is isomorphic to $\Bbb{Z}$, and the homomorphism 
 $$R_k[G] \cong \pi _0 K^G(*) \rightarrow K^G(Spec(F)) \cong \Bbb{Z}
 $$
 can clearly be identified with the augmentation $\epsilon$.  It follows that $I[G]$ annihilates $\pi _*K^G(Spec(F))$.  This result also clearly also applies to the situation where the action is on a disjoint union of schemes of the form $Spec(F)$, as one readily sees by Proposition \ref{disjoint}.  In  the general case, we have  the Gersten filtration on the category of modules over $X$, as in Theorem 5.4 of \cite{quillen}. It  is clearly invariant under the action of $G$, and further the subquotients are identified with equivariant $K$-theory spectra of disjoint unions of spectra of fields.  It follows that the homotopy groups of the  $K$-theory spectra of the subquotients in the Gersten filtration are annihilated by $I[G]$.  An easy inductive argument now gives the result.  
 \end{Proof}

We suppose we are given a scheme $X$ and a continuous action by a profinite group $G = \column{lim}{\column{\leftarrow}{i \in I}} G_i$.  Continuity here simply means that the action factors through one of the finite quotients $G_i$.  Let $I(X)$ be the full left filtering subcategory of all $i \in I$ such that the $G$-action on $X$ factors through $G_i$.  Then we define the equivariant $K$-theory spectrum $K^G(X)$ to be the homotopy colimit
$$ \column{hocolim}{\column{\rightarrow}{i \in I(X)}} K^{G_i}(X)
$$

$K^G(X)$ is a  ring spectrum for profinite $G$ as well, by the results of \cite{schwede}, and the above results clearly apply in this context as well.
\subsection{Homotopy properties}

Let $G$ be any finite group, and let $G$ act by automorphisms on a vector space $V$ over a field $k$ via a representation $\rho$.  We suppose that $\# (G)$ is relatively prime to $char(k)$. We may form the skew polynomial ring $ k[V]_{\rho} [G]$ (see \cite{goodearl}), and it is clear that the category $\mbox{Mod}^G(V)$ ($V$ is regarded as an affine scheme) is equivalent to the category of finitely generated modules over  the regular Noetherian ring $ k[V]_{\rho} [G]$.   Let $\mbox{Rep}_k[G]$ denote the category of finite dimensional $k$-linear representations of $G$, which can be identified with the category of finitely generated modules over the group ring $k[G]$.  We have an inclusion $i: k[G] \hookrightarrow k[V]_{\rho} [G]$ and a projection $\pi : k[V]_{\rho} [G] \rightarrow k[G]$ of rings.  The following result is now an immediate consequence of Theorem 7 of \cite{quillen}.

\begin{Theorem} \label{eqhomotopy} The ring homomorphisms $i$ and $\pi$ give rise to exact functors 
$$ k[V]_{\rho} [G] \column{\otimes}{k[G]} - : \mbox{\em Rep}_k[G] \rightarrow \mbox{ \em Mod}^G(V)$$ 
and 
$$k[G] \column{\otimes}{k[V]_{\rho} [G]} - :   \mbox{\em Mod}^G(V)  \rightarrow \mbox{\em Rep}_k[G]
$$
Both functors induce equivalences on $K$-theory spectra.  
\end{Theorem}  
\subsection{Completions}
We recall from \cite{completion} that given a homomorphism of commutative ring spectra $f:A \rightarrow B$, and an $A$-module $M$, one may construct the {\em derived completion} $M^{\wedge}_f$.  Here are some  examples of this construction.

{\bf Example:} Any commutative ring may be thought of as a commutative  ring spectrum via the Eilenberg-MacLane construction.  If $f: A \rightarrow B$ is a surjective homomorphism of commutative Noetherian rings, and $M$ is a finitely generated $A$-module, then the derived completion $M^{\wedge}_f$ is the Eilenberg-MacLane construction on the usual completion $M^{\wedge}_I$, where $I$ is the kernel of $f$.  

{\bf Example} Let $f: S^0 \rightarrow \Bbb{H}_p$ be the mod $p$ Hurewicz map, and let $X$ be any spectrum (and therefore an $S^0$-module).  Then $X^{\wedge}_f$ is the usual Bousfield-Kan  notion of completion of the spectrum at the prime $p$ (see \cite{bousfield}).  

{\bf Example:}  Let $A = R[\Bbb{Z}_p]$  be the representation ring of the additive group of the $p$-adic integers, defined as the direct limit of the representation rings of the finite cyclic $p$-th power order groups.  Let $f: A \rightarrow \Bbb{F}_p$ be the augmentation followed by reduction mod $p$.  Then $A^{\wedge}_f$ is equivalent to the $p$-adic completion of the integral group ring of the simplicial group of points on the circle group.  In particular, the homotopy groups are $\cong \Bbb{Z}_p$ for dimensions 0 and 1, and are $\cong 0$ otherwise.

It has a number of useful properties which can be found in \cite{completion}, and we record some of them here. Let $f: A \rightarrow B$ be  a homomorphism of commutative ring spectra.

\begin{Proposition} \label{compcofib} Suppose that 
$$  M \rightarrow N \rightarrow P
$$
is a cofibration sequence of $A$-modules.  Then the naturally defined sequence 
$$ M^{\wedge}_f \rightarrow N ^{\wedge}_f \rightarrow P ^{\wedge}_f 
$$
is also a cofibration sequence.  
\end{Proposition} 

\begin{Proposition}\label{criterion}  Suppose we are given a homomorphism $\lambda: M \rightarrow N$ of $A$-modules. Suppose further that the derived smash product map $id_B \wedge \lambda : B \column{\wedge}{A} M \rightarrow B \column{\wedge}{A}N$ is an equivalence of spectra.  Then the natural map $\lambda ^{\wedge}_f: M^{\wedge}_f \rightarrow N^{\wedge}_f$ is an equivalence of spectra.  In particular, if $ B \column{\wedge}{A}M \simeq *$, then $M^{\wedge}_{f} \simeq *$.  
\end{Proposition} 

{\bf Remark:} By the derived smash product, we mean the smash product construction applied to cofibrant replacements for $M$ and $N$.  If they are already cofibrant, then one can apply the smash product directly. 

\begin{Proposition}  Let $f: A \rightarrow B$ be a homomorphism of $(-1)$-connected commutative  ring spectra, with $\pi _0 (f) $ surjective,  and let $M$ be an $A$-module. Let $I$ denote the kernel of $\pi _0 (f)$.  Suppose that the $\pi _0A$-modules $\pi _j M$ are such that for each $j$, there is an integer $e_j$ so that $I^{e_j} \cdot \pi _j M = \{ 0 \}$.  Then the natural map $\eta : M  \rightarrow M^{\wedge}_f$ is an equivalence.
\end{Proposition}
\begin{Proof} Theorem 7.1 of \cite{completion} reduces us to the case where $A$ is  the ring $\pi _0A$, $B$ is the ring $\pi _0 B$,   and $M$ is a module for which $I^eM = \{ 0 \}$ for some $e$.  We then obtain the finite filtration 
$$ \{ 0 \} = I^e M \subseteq I^{e-1}M \subseteq \cdots \subseteq I^2M \subseteq IM \subseteq M
$$
each of whose quotients is annihilated by $I$.  It follows from (6) of Proposition 3.2 of \cite{completion} that $\eta : I^sM/I^{s+1}M \rightarrow (I^sM/I^{s+1}M )^{\wedge}_f$ is an equivalence.  A straightforward application of Proposition \ref{compcofib} now gives the result. 
\end{Proof} 

The following technical criterion will be useful in later sections. 

\begin{Proposition} \label{technical} Let 
$$ M_0 \stackrel{f_0}{\rightarrow} M_1 \stackrel{f_1}{\rightarrow}M_2 \stackrel{f_2}{\rightarrow}\cdots
$$
be a directed system of  module spectra over the commutative  ring spectrum $A$.  Suppose further that we are given a homomorphism $f : A \rightarrow B$ of commutative  ring spectra. Finally, suppose that the $Ind$-groups 
$$ \pi _t(B \column{\wedge}{A} M_0) \rightarrow   \pi _t(B \column{\wedge}{A} M_1) \rightarrow   \pi _t(B \column{\wedge}{A} M_2) \rightarrow \cdots 
$$
are $Ind$-trivial.  Then the $A$-module $M_{\infty} = \mbox{ \hspace{-.3cm}} \column{hocolim}{r} M_r$ satisfies $B \column{\wedge}{A}M_{\infty} \simeq *$, and therefore $(M_{\infty})^{\wedge}_f \simeq *$. \end{Proposition}
\begin{Proof}  Because smash products are colimit constructions, it is possible to commute the constructions $B \column{\wedge}{A} -$ and colim, which shows that $B \column{\wedge}{A} M_{\infty} \simeq *$.  
\end{Proof}

One situation in which one can  prove that the hypotheses of Proposition \ref{technical} are satisfied is treated in the following result. 

\begin{Proposition} \label{localize} Suppose $f:A \rightarrow B$ is a homomorphism of $(-1)$-connected commutative  ring spectra, which is surjective on $\pi _0$,  and suppose $L$ is a connective commutative $A$-algebra.
 Let $I_A \subseteq \pi _0 A$ denote the ideal $Ker(f)$.  
Consider now an Ind-$A$-module of the form 

\begin{equation} \label{system} 
 L \stackrel{\times \alpha _0}{\longrightarrow} L \stackrel{\times \alpha _1}{\longrightarrow} L \stackrel{\times \alpha _2}{\longrightarrow} \cdots
\end{equation} 
where $\alpha _i $ is an element of the ideal $I_A \cdot \pi _0 L$.  Then Ind-groups
\begin{equation} \label{indsystemone}
\pi _t (B \column{\wedge}{A} L )\stackrel{\pi _t(\times \alpha _0)}{\longrightarrow} \pi _t(B \column{\wedge}{A}L)\stackrel{\pi _t (\times \alpha _1)}{\longrightarrow}\pi _t (B \column{\wedge}{A}L )\stackrel{\pi _t (\times \alpha _2)}{\longrightarrow} \cdots
\end{equation}
are Ind-trivial, and therefore by   Proposition \ref{technical}, the homotopy colimit $L_{\infty} = \mbox{}\column{hocolim}{i} L$ satisfies $B \column{\wedge}{A} L_{\infty} \simeq *$.

\end{Proposition} 
\begin{Proof} From the theory of derived tensor products for modules over ring spectra, as developed in \cite{schwede}, there is a spectral sequence with $E_2$-term $Tor^{\pi _*A}_* (\pi _*B, \pi _ * L)$  converging to  $\pi _* (B \column{\wedge}{A} L)$.  We now have an inductive system of spectral sequences corresponding to diagram \ref{indsystemone} above, and  claim that maps in the system are identically zero.  To prove this, we must show that any $\lambda \in I_A \cdot \pi _0 L$ induces the zero map on all $Tor$-groups.  To see this, we write such a $\lambda $ as 
$$\lambda = \sum _s i_s \lambda _s
$$
where $i_s \in I_A$, and $\lambda _s \in \pi _0 L$.  If we can prove that each $i_s$ induces the zero map on $Tor$-groups, we will have verified the claim.  But this is clear, since multiplication by elements of $I_A$ induce the zero map on $\pi _* B$.  
 Further, the colimit of the spectral sequences converges to the colimit of the homotopy groups due to the strong convergence of the spectral sequences which is guaranteed by the connectivity hypothes on $A$ and $B$, which gives the result.   \end{Proof} 
 
 \begin{Corollary} \label{tootechnical} Suppose we are given $A,B,L, \mbox{ and } I_A,$ are as in Proposition \ref{localize},  and that $\epsilon _L : L \rightarrow B$ is a homomorphism of commutative ring spectra so that the diagram 
 
 $$
\begin{diagram}\node{A} \arrow{e} \arrow{se,b}{f}    \node{L} \arrow{s,b}{\epsilon _L} \\
\node{}   \node{B}
\end{diagram}
$$
commutes.  Let $I_L$ be the kernel of $ \pi _0 (\epsilon _L )$.  Suppose further that there is an integer $N$ so that $I_L^N \subseteq I_A \cdot \pi _0(L)$.   We suppose  we are given an Ind-$A$-module as in (\ref{system}) above, but with the assumption that $\alpha _i \in I_L$ rather than that $\alpha _i \in I_A$.  Then the Ind-groups in (\ref{indsystemone}) are trivial, and we again have that $B \column{\wedge}{A}L_{\infty} \simeq *$
 \end{Corollary} 
 \begin{Proof}  We note that any $N$-fold composite in the system (\ref{indsystemone}) has its image in $I_A \cdot \pi _0 (L)$, and a standard cofinality argument gives the result.  
 \end{Proof}

\section{Totally torsion free groups}

\noindent We say a profinite $p$-group $G$ is {\em totally torsion free}  if for every subgroup $K \subseteq G$ the quotient group $K / [K,K]$ is torsion free.  Free abelian pro-$p$ groups and free pro-$p$ groups are totally torsion free,  and free products (in the category of profinite pro-$p$ groups) of totally torsion free pro-$p$ groups are totally torsion free.  However, the condition is quite restrictive, as the following example shows.  

{\bf Example:} The {\em $2$-adic Heisenberg group}, i.e. the group of $3 \times 3$ upper triangular matrices with entries in $\Bbb{Z}_2$ and ones on the diagonal is {\em not} TTF, although it is easily seen to be torsion free.  To see this, we note that the subgroup of matrices of the form 
 
 \[  \left[ \begin{array}{ccc}
 1 & 2k  &  l  \\
 0  &  1  &  m  \\
 0  &  0  &  1  
\end{array} 
\right]
\]
has abelianization isomorphic to 

$$ \Bbb{Z}_2 \oplus \Bbb{Z} _2 \oplus \Bbb{Z}/2 \Bbb{Z}
$$
The matrix 
\[ \left[ \begin{array}{ccc}
1  &  0  &  1  \\
0  &  1  &  0  \\
0  &  0  &  1
\end{array} \right] \]
is not a commutator, but its square is. 

\begin{Proposition} \label{totorfree}
Let $k$ be a field containing all roots of unity, with separable Galois  group $G_k$.  Then $G_k$ is totally torsion free.\end{Proposition}
\begin{Proof} It clearly suffices to prove that $G_k^{ab}$ is torsion free for all fields $k$ containing all the roots of unity.  Kummer theory (see e.g. \cite{lang}, \S VI.8) tells us that $G_k^{ab}$ is the Galois group of the infinite extension $k^{ab} = \bigcup _mk((k^*)^{\frac{1}{m}})$, and therefore that any non-identity element $g \in G_k^{ab}$ acts non-trivially on an element of the form $\sqrt[m]{\kappa}$ for some $m$ and $\kappa \in k$.  Suppose now that $g^n = e$ in $G_k^{ab}$.   We have $g \cdot \sqrt[m]{\kappa} = \zeta \sqrt[m]{\kappa}$ for some $\zeta \in \mu _m (k)$.  Now choose any $n$-th root $\xi$ of $\sqrt[m]{\kappa}$.  We see that $g \cdot \xi = \eta \xi$ for some $mn$-th root of unity $\eta $ satisfying $\eta ^n = \zeta$.  Now, $g^n \cdot \xi = \eta ^n \xi = \zeta \xi \neq \xi$, which contradicts the assumption that $g$ is $n$-torsion.  
\end{Proof} 

{\bf Remark:} This structural property of Galois groups will allow us to identify completion constructions based on representation theory with $K$-spectra of certain pro-schemes which appear to be the algebraic geometric version of the classifying space construction in topology.  It would be interesting to consider how this property relates with the Bloch-Kato conjecture.

Divisibility properties of ideals in representation rings turn out to be of importance in the analysis of the conjectures which we will be discussing.  We record for use in future work one particular property of the representation rings of totally torsion free $l$-profinite groups.  First, consider the representation ring $R[\Bbb{Z}_l]$, defined as the colimit of the rings $R[\Bbb{Z}/l^n\Bbb{Z}]$.  Let $I \subseteq R[\Bbb{Z}_l]$ denote the augmentation ideal.  

\begin{Proposition}\label{div}  We have that $I/I^2$ is a $l$-divisible
group, and that
$I^k = I^{k+1}$ for all $k$. 
\end{Proposition}
\begin{Proof} We first observe that since complex representation rings of finite abelian groups are isomorphic to the group rings of their character groups, we have that $R[\Bbb{Z}_l]$ is isomorphic to the group ring $\Bbb{Z}[C_{l^{\infty}}]$, where  $C_{l^{\infty}}$ is the infinite union of the groups $\Bbb{Z}/^n \Bbb{Z}$.  It is a standard result for any group $G$ that in the
group ring $\Bbb{Z}[G]$, the augmentation ideal $I[G]$ has the property
that
$I[G]/I[G]^2 \cong G^{ab}$.  Since $C_{l^{\infty}}$ is abelian, we find
that 
$I/I^2 \cong C_{l^{\infty}}$, which is $l$-divisible.  For the second
statement, we note that $I^k/I^{k+1}$ is a surjective image of 
$\underbrace{I/I^2 \otimes \ldots \otimes I/I^2}_{k \mbox{ factors}}$. 
But it is clear that $C_{l^{\infty}} \otimes C_{l^{\infty}} = 0$, so
$I^k/I^{k+1} = 0$. This gives the result.  
\end{Proof} 

\begin{Corollary}\label{newdiv} We have that $((l)+I)^k = (l^k) + I$.
Equivalently, for any element $\theta \in I$, and any positive integers
$s$ and $t$, there is an $\eta \in I$ and $\mu \in I^k$ so that 

$$   \theta = l^s \eta + \mu . 
$$ 
\end{Corollary}

\begin{Proof} We have 
$$ ((l)+I)^k = \mathop{\Sigma} _t l^t I^{k-l} = (l^k)+l^{k-1}I +I^2. 
$$
But \ref{div} implies that $l^{k-1}I + I^2 = I$, which gives the result.  
\end{Proof}

\noindent  Recall that the representation ring of a finite group is not
just a ring, but is actually a part of a {\em Green functor}.  See \cite{Bouc} for a thorough discussion of
these objects.   A Mackey functor for a group
$G$ is a functor from a category of orbits to abelian groups.  In the
case of the representation ring, there is a functor ${\cal R}$  given on
orbits by
${\cal R}(G/K) = R[K]$.  The maps induced by projections of orbits induce
restriction maps on representation rings, and transfers induce
inductions. This functor is actually a commutative {\em Green functor},
in the sense that there is a multiplication map ${\cal R} \Box {\cal R }
\rightarrow {\cal R}$, which is associative and commutative. Also for any
finite group $G$, there is another Green functor ${\cal Z}$ given on
objects  by ${\cal Z}(G/K) =
\Bbb{Z}$, and for which projections of orbits induce the identity map and
where transfers associated to projections  induce multiplication by the
degree of the projection.  We may also consider ${\cal Z} /l {\cal Z}$,
which is obtained by composing the functor ${\cal Z}$ with the
projection $\Bbb{Z} \rightarrow  \Bbb{Z}/l\Bbb{Z}$. The {\em
augmentation} is the morphism of Green functors ${\cal R} \rightarrow
{\cal Z}$ which is given on an object $G/K$ by the augmentation of
$R[K]$. The mod-$l$ augmentation $\varepsilon$  is the composite 
${\cal R} \rightarrow {\cal Z} \rightarrow {\cal Z}/l{\cal Z}$. The
theory of Green functors is entirely parallel with the theory of rings,
modules, and ideals.  We may therefore speak of the {\em
 augmentation ideal} ${\cal I}$ and the ideal ${\cal J} = (l) + {\cal
I}$, as well as the powers of these ideals. We can therefore also speak
of completion at a Green functor ideal. We also observe that the theory
of Mackey and Green functors extends in an obvious way to profinite
groups, by considering the category of finite $G$-orbits.  T

\begin{Proposition} \label{torfreerep} Let $G$ be a totally torsion free
$l$-profinite group.  Then the natural map ${\cal R}^{\wedge}_{\cal J}
\rightarrow {\cal Z}^{\wedge}_l$ is an isomorphism. Consequently, for $G$
the maximal pro-$l$ quotient of the absolute Galois group of a geometric
field, we have that ${\cal R}^{\wedge}_{\cal J}  \cong {\cal Z}^{\wedge}_l$. 
\end{Proposition}

\begin{Proof}  We will verify that ${\cal R}^{\wedge}_{\cal J}(G/G) \rightarrow
{\cal Z}^{\wedge}_l(G/G)$ is an isomorphism.  The result at any $G/K$
for any finite index subgroup will follow by using that result for the
totally torsion free group $K$. It will suffice to show that for every
finite dimensional representation
$\rho$ of $G/N$, where
$N$ is a normal subgroup of finite index, and every choice of positive
integers
$s$ and
$t$, there are elements $x \in {\cal R}(G/G)$ and $y \in {\cal I}^t(G
/G)$ so that $[dim\rho] - [\rho] = l^s x + y$.  We recall {\em
Blichfeldt's theorem} \cite{serre}, which asserts that there is a
subgroup
$L$ of $G/N$ and a one-dimensional representation $\rho _L$ of $L$ so that
$\rho$ is isomorphic to the representation of $G/N$ induced up from $\rho
_L$. It follows that $[dim\rho] - [\rho] = i_L^{G/N}(1-\rho _L)$. Let
$\overline{L}$ denote the subgroup $\pi ^{-1}L \subseteq G$, where $\pi
: G \rightarrow G/N$ is the projection, and let $\rho _{\overline{L}} =
\rho_L \circ \pi$.   Then we clearly also have
$[dim\rho] - [\rho] = i_{\overline{L}}^{G}(1-\rho _{\overline{L}})$. 
Now, $1 - \rho _{\overline{L}}$ is in the image of $R[\overline{L}^{ab}]
\rightarrow R[\overline{L}]$, and let the corresponding
one-dimensional representation  of $\overline{L}^{ab}$ be $\rho
_{\overline{L}^{ab}}$. Since $\overline{L}^{ab}$ is abelian and torsion
free (by the totally torsion free hypothesis), we may write 
$1-\rho
_{\overline{L}^{ab}} = l^s \xi + \eta$, where $\eta \in
I^t(\overline{L}^{ab})$ and where $\xi \in R[\overline{L}^{ab}]$, by
Corollary \ref{newdiv}. This means that we may pull $\xi $ and $\eta$ back along the
homomorphism 
$\overline{L} \rightarrow \overline{L}^{ab}$, to get elements
$\overline{\xi} \in R[\overline{L}] $ and $\overline{\eta} \in
I^t(\overline{L})$ so that $\rho _{\overline{L}} = l^s \overline{\xi} +
\overline{\eta}$. Since ${\cal I}^t$ is closed under induction, we have
that $i_{\overline{L}}^G(\overline{\eta}) \in {\cal I}^t(G/G)$. Now we
have that $[dim \rho] - \rho = i_{\overline{L}}^G (1 - \rho
_{\overline{L}}) = p^s i_{\overline{L}}^G \overline{\xi} +
i_{\overline{L}}^G \eta$.  The result follows.   
\end{Proof}




\section{Actions on schemes}\label{schemeactions}

Let $k$ be our fixed separably closed  field, and let $L$ be any abelian group.  $L$ can be regarded as a discrete group scheme (see \cite{waterhouse}), and we denote this scheme by $L^0$.    We let $\Bbb{G}_m$ denote the multiplicative group scheme.  For any set $X$ and scheme $Y$, we let $F(X,Y)$ denote the scheme of functions from $X$ to $Y$.  This can be viewed as the inverse limit in the category of schemes over the discrete category with object set $X$, of the constant functor with value $Y$.  Let $G$ be any commutative group scheme.  Then we will  also write $Hom(L^0, G)$ for the subscheme of $F(L, G)$ consisting of the homomorphisms.  Of course, $Hom(-,G)$ defines a contravariant functor from abelian groups to the category of  commutative group schemes.  

We next consider automorphism groups in this context.  The group $k^*$ acts on $\Bbb{G}_m(k)$ by translation.  Consequently,  for any abelian group $L$, we have an action of the group $Hom(L,k^*)$ by automorphisms on the commutative group scheme $Hom(L^0, \Bbb{G}_m(k))$, via pointwise action.  In addition, we have the action of the automorphism group of $L$, which we denote by $GL(L)$, by the formula $(\gamma \cdot f)(l) = f(\gamma ^{-1} l)$.  These two assemble into an action of the semidirect product $GL(L) \ltimes Hom(L, k^*)$ by automorphisms on $Hom(L^0, \Bbb{G}_m(k))$, where $GL(L)$ acts in the evident way on $Hom(L, k^*)$.  The following proposition is now a straightforward verification. 

\begin{Proposition} \label{isoconstant} The homomorphism from $GL(L) \ltimes Hom(L, k^*)$ to the automorphism group  over $Spec(k)$ of the commutative group scheme  $Hom(L^0, \Bbb{G}_m(k))$ is an isomorphism when $L$ is a finitely generated free abelian group.  
\end{Proposition}

 We now also fix a prime $p$ relatively prime to the characteristic of $k$, and define a special pro-group scheme $\frak{G}_m(k)$  on objects by setting ${\frak{G}}_m(k)(n) = \Bbb{G}_m(k)$, and on morphisms by setting ${\frak{G}}_m(k)(n_1 \geq n_2)$ equal to the $p^{n_1 -n_2}$-th power map on  $\Bbb{G}_m(k)$.   We write $Aut({\frak{G}}_m)$ for the automorphism group of $\frak{G}_m$, and we note that it is the inverse limit of a pro-group, defined as follows.  Let $A_n \subseteq Aut ({\frak{G}}_m)$ be the subgroup of elements $\{ \alpha _i \}_i$ such that $\alpha _i = Id_{\Bbb{G}_m }$ for $i \leq n$.  Each $A_n$ is a normal subgroup of $Aut({\frak{G}}_m)$, and $Aut ({\frak{G}}_m)$ is identified with the inverse limit of the pro-group 

$$  \cdots \rightarrow Aut({\frak{G}}_m)/A_{n+1} \rightarrow Aut({\frak{G}}_m)/A_{n} \rightarrow Aut({\frak{G}}_m)/A_{n-1} \rightarrow \cdots 
$$

 We define a group $U_k$ to be the inverse limit of the group of special pro-abelian group of $k$-valued points of $\frak{G}_m(k)$ and let $T_k \subseteq U_k$ be the inverse limit of the special pro-abelian group of $k$-valued points of finite $p$-th power order within $U_k$, so we have  

$$ T_k =  \cdots \stackrel{\times p}{\rightarrow} \mu _{p^{\infty}} \stackrel{\times p}{\rightarrow} \mu _{p^{\infty}}   \stackrel{\times p} {\rightarrow} \mu _{p^{\infty}} 
$$
  We observe that there is a homomorphism $i: U_k \rightarrow Aut({\frak{G}}_m)$, which sends an element of $U_k$, which can be viewed as a vector of elements $\{ \kappa _i \}$ in $k^*$ so that $\kappa _i ^p = \kappa _{i-1}$ for all $i$, to the automorphism which is given by multiplication by $\kappa _n$ on ${\frak{G}}_m(n)$.  We now obtain the following  counterpart to Proposition \ref{isoconstant}. 

\begin{Proposition} \label{repdef}Given any abelian group $L$, we may construct the diagram of group schemes $Hom(L^0, {\frak{G}}_m(k))$.  There is a natural homomorphism from the group $Hom(L,U_k)$ to the automorphism group of $Hom(L^0, {\frak{G}}_m(k))$, and as before, $GL(L)$ also acts by automorphisms on this group.  The two fit together to provide a homomorphism $\theta: GL(L) \ltimes Hom(L, U_k) \rightarrow Aut(Hom(L^0,{\frak{G}}_m (k))$, which has the following properties. 
\begin{enumerate}
 \item{$\theta$  is a continuous  isomorphism of pro-groups when $L$ is finitely generated and free.}
 \item{The subgroup $im(U_k) \subseteq Aut({\frak{G}}_m(k))$ can be characterized as the subgroup which induces the identity on the quotient group $A[Hom(L^0, {\frak{G}}_m(k))]^*/k*$, where $A[-]$ denotes the affine coordinate ring functor, and $(-)^*$ denotes the units.  We will refer to it as the {\em translation subgroup}.  Similarly, we will call the subgroup which induces a permutation action on $A[Hom(L^0, {\frak{G}}_m(k))]^*/k* $ the {\em monomial subgroup}.  Notice that this group can be characterized as the subgroup which preserves a choice of free generating monoid for the group of units in the affine coordinate ring, and consequently extends in a natural way to an affine space.   } \end{enumerate} 
\end{Proposition} 

We note that $GL(L) \ltimes Hom(L, U_k)$ contains the subgroup $GL(L) \ltimes Hom(L, T_k)$, since $T_k$ is invariant under the $GL(L)$ action on $Hom(L,U_k)$.  
\begin{Definition}\label{creps}
Let $\Gamma$ denote a pro-$p$ group, that is a  pro-object in the category of finite $p$-groups.  Then by a {\em ${\frak{G}}$-representation}, we will mean a continuous homomorphism from $\Gamma $ to the pro-group $Aut(Hom(L^0, {\frak{G}}_m(k)))$, for some finitely generated free abelian group $L$.  Of course,
$Hom(L^0, {\frak{G}}_m(k)) $ is abstractly isomorphic to ${\frak{G}}_m(k) ^n$, where $n$ is the rank of $L$.  For any $\frak{G}$-representation $\rho$, we will denote by $\rho[k]$ the representation at level $k$, i.e. the action of the $\Gamma$ on the $k$-th entry in the tower.  
\end{Definition} 

We wish to define certain one-dimensional ${\frak{G}}$-representations of $\Gamma$.  We suppose that we are given a family of continuous $k^*$-valued one-dimensional $k$-valued characters $\{ \chi _n \}$, satisfying the condition $\chi _n^p  = \chi _{n-1}$.  The group $k^*$ is equipped with the discrete topology, so continuity of a character simply means that   it factors though one of the finite quotients defining $\Gamma$.  Together, they yield a homomorphism from $\Gamma$ into $U_k$, and due to the continuity the homomorphism actually factors through $T_k$.  In summary, a sequence of continuous characters $\{ \chi _n \}$ of $\Gamma$ satisfying the compatibility condition $\chi _n^p = \chi _{n-1}$ defines a continuous  action of $\Gamma$ on the diagram ${\frak{G}}_m(k)$, i.e. a continuous homomorphism from $\Gamma$ to the pro-group $Aut({\frak{G}}_m)$.  It will be called a {\em translation representation}   since it factors through the translation subgroup.  Similarly, we will call an action which factors through the monomial subgroup a {\em monomial representation}. 

\begin{Definition} Let $k$ denote a separably closed field.  Then by a {\em $\chi$-sequence in $k^*$}, we will mean a sequence  $\{ \mu _i \}_{i \geq 1}$ of elements of $T_k$ so that $\mu ^p_i = \mu _{i-1}$.  Associated to a $\chi$-sequence $\{ \mu _i \}$, we define the one dimensional  ${\frak{G}}$-representation $\rho ( \{ \mu _i \} )$ as above.  
\end{Definition}

We now have the following definition. 

\begin{Definition} Let $G$ denote any profinite group, and suppose that we are given a continuous action of $G$ on a scheme $X$ , where the automorphism group of $X$ is regarded as discrete.  Letting $K \subseteq G$ denote any closed subgroup, we say the action is {\em $K$-principal} if the stabilizer of $G_x = K$ for all closed points $x \in X$, where $G_x$ denotes the stabilizer group of $x$.  \end{Definition}

We now suppose we have a continuous action of a profinite group on the  special pro-group scheme $Hom(L^0, {\frak{G}}_m(k))$, say $\rho$.  This means that we have a family of homomorphisms $$\rho _n : G \rightarrow Aut(Hom(L^0, \Bbb{G}_m(k))$$ so that the diagrams 

$$
	\begin{diagram}
	\node{ \Bbb{G}_m(k)} \arrow{s,b} {\times p^{n_1 - n_2}}\arrow{e,t}{\rho  _{n_1}(g)} \node{ \Bbb{G}_m(k)}\arrow{s,b}{\times p^{n_1 - n_2}}\  \\
	\node{  \Bbb{G}_m(k)} \arrow{e,t}{\rho  _{n_2}(g)} \node{ \Bbb{G}_m(k)}
	\end{diagram} 
	$$
commute for all $g \in G$ and all $n_1 \geq n_2$.   Let $K \subseteq G$ be any closed subgroup.  Then we say the action $\rho$ is {\em $K$-principal } if each action $\rho _n$ is $K \cdot Ker(\rho _n)$-principal.  We record a useful property of $K$-principal actions.  

\begin{Proposition}\label{principalcriterion}  Let $\rho _1$ and $\rho _2$ be continuous actions of a profinite group on  $Hom(L^0_1, \frak{G}_m(k))$ and $Hom(L^0_2, \frak{G}_m(k))$, respectively, where $L_1$ and $L_2$ are finitely generated abelian groups.  Then we have the evident product action $\rho _1 \times \rho _2$ on $Hom (L_1^0 \times L_2 ^0, \frak{G}_m(k)) \cong Hom(L^0_1, \frak{G}_m(k)) \times Hom(L^0_2, \frak{G}_m(k))$.  Let  $\{ e \} \subseteq K^{\prime} \subseteq K \subseteq G$ be closed subgroups.  Suppose that $\rho _1 $ is $K$-principal, and that $\rho _2 | K$ is $K^{\prime}$-principal.  Then $\rho _1 \times \rho _2$ is $K^{\prime}$-principal.  \end{Proposition} 
\begin{Proof} Clear.  
\end{Proof}

Next, we must consider the process of induction of actions $\rho$  on ${\frak{G}}_m(k)$ by subgroups $K$ of finite index in $\Gamma$ to actions of $\Gamma$ on finite products of ${\frak{G}}_m(k)$.  Let $G$ be a finite group, and let $K \subseteq G$ be a subgroup.  Suppose further that $X$ is any scheme equipped with a left $K$-action.  We can now consider the scheme $F(G, X)$, which is equipped with a left $K$-action defined by $(\kappa \cdot f)(g) = \kappa f(\kappa ^{-1} g)$.  Its fixed point scheme is the scheme of $K$-equivariant maps from $G$ to $X$, and we denote it by $F^K(G,X)$.  $F^K(G,X) $ is now equipped with a left $G$-action, defined by $(\gamma \cdot f) (g)   = f(g \gamma )$, and we write $i_K^G (X)$ for it. It is clear from the definition that this construction extends to ${\frak{G}}$-representations, in the sense that given a profinite group $G$, a closed subgroup of finite index $K$, and a ${\frak{G}}$-representation $\rho$ of $K$, there is a natural $G$-representation $i_K^G (\rho) $ of $G$.  Its salient properties are given in the following result. 

\begin{Proposition}  \label{induction} Let $G$ be a  profinite group, and let $K \lhd G$ be a normal subgroup of finite index.  Suppose we are given a one-dimensional ${\frak{G}}$-representation $\rho$ of $K$, and construct the induced ${\frak{G}}$-representation of $G$ $i_K^G(\rho )$.  Then we have the following. 

\begin{itemize}
\item{The restriction of $i_K^G(\rho)$ to $K$ is a translation representation. }
\item{The restriction of $i_K^G(\rho)$ to $K$ is $K^{\prime}$ principal, where $K^{\prime}$ is defined as follows.  For each element $\gamma \in G/K$, we denote by $\rho ^{\gamma}$ the result of composing $\rho$ with the conjugation action of a coset representative of $\gamma$.  Then $K^{\prime}$ is given by 

$$ K^{\prime} = \bigcap _{\gamma  \in G/K} Ker(\rho ^{\gamma} )
$$
}
\end{itemize} 
\end{Proposition} 

\section{Universal pro-schemes } \label{universal} 

Let $G$ be any $p$-profinite group.  The monomial $\frak{G}$-representations of $G$ are given as continuous homomorphisms $\rho : G \rightarrow \Sigma _n \ltimes \Bbb{Z}_p^n$.  We let ${\cal R}(G) $ denote the set 
$$    {\cal R}(G) = \coprod_n {\cal R}_n(G)  = \coprod _{n} Hom^{\mbox{c}}(G, \Sigma _n \ltimes \Bbb{Z}_p^n)
$$
where $Hom ^{\mbox{c}}$ denotes the set of continuous homomorphisms.   We consider the group scheme

$$\Gamma (G) =  \prod_n \prod_{\rho \in {\cal R}_n(G)} F(\Bbb{N}, \Bbb{G}_m(k)^n) =  \prod_n \prod_{\rho \in {\cal R}_n(G)}(\Bbb{G}_m(k)^n_{\rho})^{\infty}
$$

$\Gamma (G)$ is equipped with a $G$ action where the action on $(\Bbb{G}_m(k)^n_{\rho})^{\infty}$ is the diagonal action of infinitely many copies of $\rho$.   

For any finite subset $S = \{ (\rho _1, k_1 ), \ldots , (\rho _s , k_s ) \}$ of the set ${\cal R}(G) \times \Bbb{N}$, we  write $\Gamma (G;S) $ for the summand 
$$ \prod_{i} \Bbb{G}_m(k)^{n_i}_{k_i}
$$
where $n_i$ denotes the dimension of $\rho _i$, and where $ \Bbb{G}_m(k)^{n_i}_{k_i}$ is the copy of $\Bbb{G}_m(k)^{n_i}$ corresponding to $k_i$.  $\Gamma (G;S)$ is a finite product of copies of $\Bbb{G}_m(k)$. Further, $\Gamma (G;S)$ is closed under the action of $G$.   We have evident projections  $\Gamma (G;S^{\prime} ) \rightarrow \Gamma (G;S) $ whenever $ S \subseteq S^{\prime}$, which are $G$-equivariant.  Let ${\cal F} = {\cal F}({\cal R}(G) \times \Bbb{N})$ denote the partially ordered set  of finite subsets of 
${\cal R}(G) \times \Bbb{N}$.  ${\cal F}^{op}$ is now a left filtering partially ordered set, and  therefore the functor $S \rightarrow \Gamma (G;S)$ is now a pro-$G$-scheme, which we will denote by ${\cal E}G$.  For any pro-$G$-scheme ${\cal S}$, we define  $K^G({\cal S})$ as the corresponding homotopy colimit of $K$-theory spectra.   Given any pro-scheme ${\cal S}$, we can define the product pro-$G$-scheme ${\cal E}G \times {\cal S}$ in the usual way.   

\begin{Proposition} Let ${\cal G}(\rho)$ denote the pro-$G$-scheme associated to any monomial $\frak{G}$-representation $\rho$.  Then the natural map 
$$  K^G({\cal E}G) \rightarrow K^G({\cal E}G \times {\cal G}(\rho))
$$
is an equivalence of spectra.  
\end{Proposition} 
\begin{Proof}  Immediate from the definition of ${\cal E}G$.  
\end{Proof}

We now need  a technical lemma. 

\begin{Lemma} \label{principalfinal} Let $G$ be a totally torsion free pro-$p$-group, and let ${\cal F}$ be as above.  Let ${\cal F}^{prin}\subseteq {\cal F}$ be the partially ordered subset consisting of those $S \in {\cal F}$ for which $\Gamma (G;S)$ is $K$-principal for some closed subgroup $K \subseteq G$.  The ${\cal F}^{prin}$ is final in ${\cal F}$.  
\end{Lemma} 
\begin{Proof} We must prove that for any $S \in {\cal F}$, there is a $S^{\prime} \in {\cal F}^{prin}$ with $S \subseteq S^{\prime} $.  First, we observe that for any $\frak{G}$-representation $\rho: G \rightarrow \Sigma _n \ltimes \Bbb{Z}_p^n$, if the projection $$G \stackrel{\rho}{\rightarrow} \Sigma _n \ltimes \Bbb{Z}_p^n \rightarrow \Sigma _n$$
is trivial, then $\rho $ is $K$-principal, where $K$ is the kernel of $\rho$.  This is immediate since the action is by translation in the commutative algebraic group.  
We next observe that it will suffice to prove that for every closed normal subgroup $N \lhd G$ of finite index, there is a principal $\frak{G}$-representation of $G$ which is $K$-principal for some closed normal subgroup   $K\subseteq N$.  For, consider any $\frak{G}$-representation $\rho : G \rightarrow \Sigma _n \ltimes \Bbb{Z}_p^n$.    We let $N_{\rho}$ denote the kernel of $\pi \compcirc \rho$, which is a  closed normal subgroup of finite index.  Let $N_{\rho}^{\prime}$ be the kernel of $\rho $ itself.   If $\sigma $ is a $K$-principal $\frak{G}$-representation with $ K \subseteq  N_{\rho}$, then we claim $\rho \times \sigma$ is a $K \cap N^{\prime}_{\rho}$-principal $\frak{G}$-representation.  To see this, we note that  $\rho |N_{\rho}$ is $N_{\rho}^{\prime} $-principal from the first statement. Further, the stabilizer of any point in $\rho \times \sigma $ must be contained in $N_{\rho}$ since $K \subseteq N_{\rho}$.  The conclusion is that the stabilizer is $K \cap N_{\rho}^{\prime}  $, which is what was required.  

To construct a  $\frak{G}$-representation which is $K$-principal for some $K \subseteq N$,  we  assume by induction that we have constructed such $\frak{G}$-representations for all subgroups of index less than the index of $N$.  We consider the finite list $\{ N_i \}_{i= 1}^s$ of all normal subgroups of $G$ which contain $N$ properly, and for each one we construct a $K_i$ principal $\frak{G}$-representation $\rho _i$ of $G$, where $K_i \subseteq N_i$.   The $\frak{G}$-representation $\rho _1 \times \rho _2 \times \cdots \rho _s$ is now $K_1 \cap \cdots \cap K_s$-principal.  We now consider the subgroup $N \cdot (K_1 \cap  \cdots \cap  K_s) $.   This subgroup contains $N$, and therefore corresponds to a subgroup of the quotient $G/N$, contained in $N_1/N \cap N_2/N \cap \cdots \cap N_s/N$.  This group can be characterized as the intersection $\gamma (G/N)$ of all non-trivial normal subgroups in $G/N$.  Since $N$ is a closed normal subgroup of the pro-$p$ group $G$, we know that $G/N$ is a $p$-group, and hence its center $Z(G/N)$ is non-trivial.  Hence, $\sigma (G/N) \subseteq Z(G/N)$ is a central subgroup, which implies that every subgroup of $\sigma (G/N)$ is normal in $G/N$.  By the definition of $\sigma (G/N)$, this means that $\sigma (G/N)$ has no non-trivial proper subgroups, and the only such $p$-groups are $\{ e \}$ and $\Bbb{Z}/p \Bbb{Z}$. 
 Consequently, either 
$N \cdot (K_1 \cap  \cdots \cap  K_s) /N$ is trivial, in which case $K_1 \cap  \cdots \cap  K_s \subseteq N$, and we are done, or it is cyclic of order $p$.  In this case, we are given a homomorphism from $ \overline{N} = N \cdot(K_1 \cap \cdots \cap K_s)$ to $\Bbb{Z}/p\Bbb{Z}$, and therefore a character $\chi _1$ of   
$\overline{N}$.  Since the group $G$ is totally torsion free, we find that the $\overline{N}^{ab}$ is torsion free, and therefore that the character group is divisible.  Construct a sequence of characters $\chi _i$ so that $\chi _i ^p = \chi _{i-1}$, with $\chi _1$ the previously defined character.  This sequence constructs a one-dimensional $\frak{G}$-representation  of $\overline{N}$, which we now induce up to a $\frak{G}$-representation $\overline{\rho}$ of dimension equal to the index of $\overline{N}$ in $G$.  We claim that the $\frak{G}$-representation $\overline{\rho} \times \rho _1 \times \cdots \times \rho _s$ is $Ker(\overline{\rho}) \cap (K_1 \cap K_2 \cap \cdots \cap K_s)$-principal.  This follows immediately from Proposition \ref{principalcriterion}   because the stabilizer of any point in  $\overline{\rho} \times \rho _1 \times \cdots \times \rho _s$ is contained in $\overline{N}$, and therefore in $\overline{N}$ and $Ker(\overline{\rho})$. On the other hand, it is clear that any element of  $Ker(\overline{\rho}) \cap (K_1 \cap K_2 \cap \cdots \cap K_s)$ acts trivially, which gives the result.  \end{Proof} 

Let ${\cal S}$ denote the partially ordered set ${\cal F} \times \Bbb{N}$, with the partial order on $\Bbb{N}$ being the usual one, viewed as a category with a unique morphism $m \rightarrow n$ whenever $m \geq n$, and no morphisms from $m$ to $n$ otherwise.  For any element $S = (\{ (\rho _1, m_1) , \ldots, (\rho _s, m_s ) \}, n)$, we define a pro-$G$-group scheme $\Psi$ by setting 
$$ \Psi (S) = \prod _{i=1}^s F( \Bbb{N} , \Bbb{G}_m(\rho _i [n]))_{m_i}
$$
where we recall that for any $\frak{G}$-representation $\rho$, $\rho [k]$ denotes the value of $\rho$ at level $k$, and where $\Bbb{G}_m(\rho[k])$ denotes the group scheme on which $\rho [k]$ acts.  The maps in the pro-object are evident.  We can apply $K^G$ to this pro-$G$-scheme, to obtain an Ind-object $K^G(\Psi)$, and it is clear from the definition that $K^G(\Psi ) \cong K^G({\cal E}G)$, since $K^G$ applied to any $\frak{G}$-representation $\rho$  is defined to be the colimit of $K^G(\rho [k])$ over $k$.  
We define  ${\cal B}G  = {\cal E}G /G$ to be the pro-scheme obtained by performing the orbit space construction objectwise.  In particular, the orbit pro-scheme of a $\frak{G}$-representation $\rho$ is the pro-scheme 
 $k \rightarrow \rho [k]/G$.  There is of course a natural map ${\cal E}G \rightarrow {\cal B}G$, and a corresponding map of $K$-theory spectra $K^G({\cal B}G) \rightarrow K^G({\cal E}G)$.  If we precompose this with the map $K({\cal B}G) \rightarrow K^G ({\cal B}G)$ obtained by considering the inclusion of the modules with trivial action into all $G$-modules, we have a natural map 
 $$ \eta : K({\cal B}G) \rightarrow K^G({\cal E}G)
 $$

\begin{Proposition}\label{orbit}  Let $W$ be any smooth Noetherian scheme with continuous $G$-action.  Then the map $K({\cal E}G \column{\times}{G}  W) \rightarrow K^G ({\cal E}G \times W)$ is an equivalence of spectra.  In particular, the natural map $K({\cal B}G) \rightarrow K^G({\cal E}G)$ is an equivalence.  
\end{Proposition} 
\begin{Proof}  We prove the case $W = *$.  The extension to more general $W$  is routine.  For any $S  = (\{ (\rho _1, m_1) , \ldots, (\rho _s, m_s ) \}, n) \in {\cal S}$, we let $N(S) \subseteq G$ be the maximal normal subgroup of $G$ so that the action $\rho _i[n]$ factors through $G/N(S)$ for all $i = 1, \ldots , s $.  Let ${\cal S}^{\prime} \subseteq {\cal S}$ denote the subset of all $(\{ (\rho _1, m_1) , \ldots, (\rho _s, m_s ) \}, n)$ so that there is an $i$ so that the action $\rho _i$ is $N(S)$-principal.  It follows from Lemma \ref{principalfinal} that ${\cal S}^{\prime}$ is final in ${\cal S}$, and therefore that the homotopy colimit which computes $K^G({\cal E}G  )$ may be computed over ${\cal S}^{\prime}$ instead.  There is also a new  spectrum valued functor  functor ${\cal K}^G$ defined on   ${\cal S}^{\prime}$, given by $S \rightarrow K^{G/N(S)}(\Gamma(G;S)[n])$, and a natural transformation ${\cal K}^G  \rightarrow K^G( \Gamma(G;S))$, obtained by applying contravariant  functoriality of $K^G$ in $G$ to the homomorphism $G \rightarrow G/N(S)$, together with the homomorphism $\rho \rightarrow \rho [n]$.   There  is a third functor defined on ${\cal S}^{\prime}$, which we denote by ${\cal O}$, and which is given by ${\cal O} = K(\Gamma (G;S)/G)$.  It is equipped with a natural transformation (as with $\eta$ above) ${\cal O} \rightarrow {\cal K}^G$.  The composite ${\cal O} \rightarrow K^G(\Gamma(G;S))$ can be seen to induce an equivalence on homotopy colimit spectra using the finality of ${\cal S}^{\prime}$ in ${\cal S} $ together with Proposition \ref{free}.    The same finality shows readily that $\column{hocolim}{\rightarrow }{\cal O}$ is equivalent to $K({\cal B}G)$, which gives the result.  
\end{Proof} 

Finally, we will prove a technical result which will be useful in analyzing completions of $K^G ({\cal E}G)$. 

\begin{Proposition} Let $G$ be a totally torsion free $p$-profinite group.  Then the action of $\pi _0 K^G(*) \cong R[G]$ satisfies $I[G] \cdot \pi _* K^G({\cal E}G) = \{ 0 \}$.  
\end{Proposition} 
\begin{Proof}  We first observe that the natural ring homomorphism  $j: R[G]\rightarrow \pi _0 K^G({\cal E}G)$ is surjective.  In view of the colimit definition of $K^G({\cal E}G)$, it will suffice to prove that given any homomorphism $ \rho: G \rightarrow \Sigma _n \ltimes \Bbb{Z}^n_p \hookrightarrow Aut(\Bbb{G}_m(k)^n)$, the map $R[G] \rightarrow K^G(\Bbb{G}_m(k)^n _{\rho})$  is surjective.  But from the homotopy property Theorem \ref{eqhomotopy}, it follows that this map is induced by the map of schemes with $G$-action $\Bbb{G}_m(k)^n _{\rho} \rightarrow \Bbb{A}_m(k)^n _{\rho}$, where $\Bbb{A}_m(k)^n _{\rho}$ denotes the extension of the $G$-action $\rho$ to the affine space. This means that the map is the $\pi _0$-part of a localization sequence for $K^G$-theory, and is therefore surjective. 
If we can show that $j(I[G]) = \{ 0 \}$, we will be done.  This means that we will need to show that for any representation $\rho $ of $G$, we have that $[\rho ] - dim(\rho)$ vanishes in $K^G({\cal E}G)$.  Blichfeldt's theorem (see \cite{serre}) tells us that  any representation of $G$ is monomial, i.e it factors through the subgroup $\Sigma _n \ltimes (k^*)^n \hookrightarrow GL_n(k)$.   We suppose we are given a representation $\rho$  and such a factorization of it.  Associated to the factorization, we have the associated permutation representation $\overline{\rho}$  obtained by composing the factorization with the projection $\Sigma \ltimes (k^*)^n \rightarrow \Sigma _n$.  We claim that the element $[\rho] - [\overline{\rho }]$ vanishes in $K^G({\cal E}G)$.   Let $\rho$ denote any  monomial representation, and construct the corresponding actions on an affine space $\Bbb{A}^n(k) _{\rho}$ and on a  product of $\Bbb{G}_m(k)$'s $\Bbb{G}_m(k)^n _{\rho}$.  We consider the portion 

$$ \pi _0 K(\mbox{Mod}^G(\Bbb{A}^n(k) _{\rho})^{{tor}}) \stackrel{\lambda}{\rightarrow }
\pi _0 K^G(\Bbb{A}^n(k) _{\rho}) \rightarrow \pi _0 K^G(\Bbb{G}_m(k)^n _{\rho}) \rightarrow 0
$$
of the localization sequence, where $\mbox{Mod}^G(\Bbb{A}^n(k) _{\rho})^{{tor}}$ denotes the category of all $\Bbb{A}^n(k) _{\rho}$-modules with $G$-action compatible with the $G$-action on $\Bbb{A}^n(k) _{\rho}$, and which are $x_1 x_2 \cdots  x_n$-torsion, where the $x_i$'s are the coordinates on $\Bbb{A}^n(k) _{\rho}$ in which the representation $\rho$ is monomial.  The homotopy property shows that $R[G]  \cong \pi _0 K^G(\Bbb{A}^n(k) _{\rho})$, so it suffices to show that $[\overline{\rho}] - [\rho]$ is in the image  of $\lambda$.  We note that in proving this fact, we may now deal with the finite quotient of $G$ which is the image of $\rho$, so that all constructions can assume that $G$ is a  finite $p$-group, which we now do.   In order to prove the claim, we need to make a definition. First, we let $A$ be the affine coordinate ring of $\Bbb{A}^n(k) _{\rho}$.  It is the symmetric algebra on the representation $\rho$, and as such is equipped with a $G$-action.   Let $H \subseteq G$ be any subgroup.  Let $M$ be any $H$-twisted $A$-module, and let $A[G]_{\rho}$ and $A[H]_{\rho}$ be the twisted $A$-group rings corresponding to the representation $\rho$.    Then by the induced module $i_H^G(M)$, we will mean the $A[G]_{\rho}$-module $A[G]_{\rho}  \column{\otimes}{A[H]_{\rho} } M$.  Since $\rho$ is monomial, there is a corresponding homomorphism $\overline{\rho} : G \rightarrow \Sigma _n$.  Let $\Sigma _n(1)$ denote the stabilizer of the element $1 \in \{ 1, \ldots , n \}$, and let $H = \overline{\rho}^{-1} (\Sigma _n (1))$.  Then the module $ M = A/x_1 A$ is an $H$-twisted $A$-module.  It follows that we can construct the $G$-twisted $A$ module $i_H^G(M)$.   On the other hand, we have the resolution 
$$ 0 \rightarrow A ^{\prime} \stackrel{x_1 \cdot}{\rightarrow} A \rightarrow M \rightarrow 0
$$
of $H$-twisted $A$-modules,  where $A$ denotes the ring $A$ with its usual $H$-action, and where $A^{\prime}$ denotes the cyclic $A$-module with the $H$ action afforded by the one-dimensional representation of $H$ on the subspace spanned by $x_1$.  We may apply $i_H^G$ to it, to obtain a resolution
$$  0 \rightarrow i_H^G (A^{\prime})  \stackrel{ i_H^G(x_1) \cdot}{\rightarrow}  i_H^G(A ) \rightarrow  i_H^G(M )\rightarrow 0
$$
of $G$-twisted $A$-modules.  It follows by Quillen's additivity theorem (Corollary 1 of Section 3 of \cite{quillen}) that the element 
$[i_H^G(A ) ] - [ i_H^G (A^{\prime}) ]$ maps to zero in $K^G_0(\Bbb{G}_m(k)^n_{\rho})$.  On the other hand,  the homotopy property shows that this element is equal to the element $[i_H^G(k\column{\otimes}{A} A)] - [i_H^G(k\column{\otimes}{A} A^{\prime} )] $, where in this case $i_H^G$ can be identified with ordinary induction of representations.  This element is now easily seen to be $[\overline{\rho}] - [\rho] \in R[G]$.  Note that the special case $n=1$ of this result shows that one-dimensional representations $\chi$ satisfy the condition that  $ [\chi] - dim(\chi)$ vanishes in $K^G ({\cal E}G)$.  To see that this result implies that $[\rho] - dim(\rho)$ vanishes in $K^G({\cal E}G)$, we assume that $\rho $ is a representation of minimal dimension $n$ (which must be $>1$)  for which $[\rho] - [\overline{\rho}] \neq 0$ in $K^G ({\cal E}G)$.   Consider the representation $\overline{\rho} $.  It is the permutation representation of $G$ on the cosets $G/H$, where as above $H = \overline{\rho}^{-1}(\Sigma _n(1))$.  The representation $\overline{\rho }$ is therefore decomposable, since it contains the one-dimensional invariant subspace spanned by the element 
$\sum _{x \in G/H} x$.  It follows that all the indecomposable summands $\{ \overline{\rho}_1, \ldots, \overline{\rho}_s \}$ of $\overline{\rho}$ have dimension strictly $< n$, which means that $[\overline{\rho }_i ] - dim(\overline{\rho }_i)$ vanishes in $K^G{\cal E}G)$ by the hypothesis on $\rho$.  But, this means that 
$$ \sum _i [\overline{\rho }_i ] - dim(\overline{\rho }_i) = [\overline{\rho}] - dim (\overline{\rho})
$$
vanishes in $K^G({\cal E}G)$, which implies that $[\rho] - dim(\rho)$ vanishes in $K^G({\cal E}G)$, which contradicts the choice of $\rho$, and therefore gives the result.  
 \end{Proof} 
 
 Now we consider the diagram of commutative  ring spectra
 $$ \begin{diagram}
 \node{S^0} \arrow{e} \arrow{s,t}{r}   \node{K^G(*)} \arrow{s,b}{\epsilon}  \\
 \node{\Bbb{H}_p} \arrow{e,t} {=}\node{\Bbb{H}_p}
 \end{diagram} 
 $$
 where $r$ denotes the ring map induced by the fundamental class in the sphere spectrum, and $\epsilon$ is the augmentation.  
 The functorality of the derived completion construction means that the diagram induces a map 
 $$\beta : K^G({\cal E}G) ^{\wedge}_p = K^G({\cal E}G) ^{\wedge}_r  \rightarrow K^G({\cal E}G) ^{\wedge}_{\epsilon} $$
\begin{Corollary} \label{nonequiv} The map $\beta $ is an equivalence of spectra.  
\end{Corollary} 
\begin{Proof}
 The algebraic-to-geometric spectral sequence (Theorem 7.1 of \cite{completion}) and the preceding Proposition show that it suffices to prove that the map on derived completions $\hat{\beta}: M^{\wedge}_p \rightarrow M^{\wedge}_{\epsilon}$ induced by the diagram 

 $$ \begin{diagram}
 \node{\Bbb{Z} } \arrow{e} \arrow{s,t}{r}   \node{R[G]} \arrow{s,b}{\epsilon}  \\
 \node{\Bbb{F}_p} \arrow{e,t} {=}\node{\Bbb{F}_p}
 \end{diagram} 
 $$
 is an equivalence for modules $M$ with trivial action by $I[G] \subseteq R[G]$, or equivalently which are obtained by restricting scalars along the homomorphism $\epsilon : R[G] \rightarrow \Bbb{Z}$.  In order to prove this, we regard the ring $R[G]$ as a commutative  ring spectrum, and let $\tilde{\Bbb{Z}}$ denote a cofibrant replacement for the 
 $R[G]$-algebra $\Bbb{Z}$, with algebra structure given by $\epsilon$, in the model category of commutative  ring spectra discussed in \cite{shipley}.  We have the natural structure map $R[G] \rightarrow \tilde{\Bbb{Z}}$, as well as a homomorphism $\tilde{\Bbb{Z}} \rightarrow R[G]$ of commutative  ring spectra obtained by extending the homomorphism of rings $\Bbb{Z} \rightarrow R[G]$ to the cofibrant replacement $\tilde{\Bbb{Z}}$.  We also have a  homomorphism $\tilde{\Bbb{Z}} \rightarrow \Bbb{H}_p$ which is compatible with the mod-$p$ reduction of the augmentation on $R[G]$.    We now have two monads $S$ and $T$ on the category of $\tilde{\Bbb{Z}}$-module spectra, where $S(M) = \Bbb{H}_p \column{\wedge}{\tilde{\Bbb{Z}}} M$ and $T(M) = \Bbb{H}_p \column{\wedge}{R[G]} M$.  Further, there are natural transformations $T \rightarrow S$ and $S \rightarrow T$ induced by the ring homomorphisms $ \tilde{\Bbb{Z}} \rightarrow R[G] \rightarrow \tilde{\Bbb{Z}}$.  It follows that the natural transformations of triples induce equivalence on the total spectra of the cosimplicial  objects of the triples $S$ and $T$, by Theorem 2.15 of \cite{completion}.  These total spectra attached to $S$ and $T$ are the derived completions over $\tilde{\Bbb{Z}}$  and $R[G]$ respectively.   The fact that weak equivalences of commutative ring spectra induce equivalences of derived completions now gives the result.  
 \end{Proof}






\section{Algebraic properties of certain representation rings}

In this section, we analyze some particular properties of the representation rings of certain semidirect products of  symmetric groups with profinite groups of the form $\Bbb{Z}_p^n$.  The results of this section are the key technical tools which allow us to prove our main results.  

\begin{Theorem} \label{general} Let $k$ be a field of characteristic $p$, and let $\epsilon: A \rightarrow k$ be an augmented $k$-algebra with augmentation ideal $I_A$. Let $G \subseteq G^{\prime} \subseteq \Sigma _n$ denote subgroups of $p$-power order.   Let $ \displaystyle B = \mathop{\otimes} ^{n}_{k} A$, equipped with the action of $\Sigma _n$, and therefore of $G$ and $G^{\prime}$. The algebra $B$ is given the tensor product augmentation $\epsilon _B$, and   the algebras $B^G$ and $B^{G^{\prime}}$ both obtain augmentations by restriction of $\epsilon_B$, and we denote the restrictions by $\epsilon _{B^G}$ and $\epsilon _{B^{G^{\prime}}}$ respectively.  Further, we let $I_{B^G}$ and $I_{B^{G^{\prime}}}$ denote the corresponding augmentation ideals.  Finally, suppose that every element in $A$ has a $p$-th root.  Then the equality of ideals 
$$ I_{B^G} = I_{B^{G^{\prime}}}\cdot B^{G}
$$
holds. 
\end{Theorem}
\begin{Proof} We first prove that the $k$-algebra $B^{G}$ also has the property that every element admits a $p$-th root.  To see this, we let $\frak{A}$ denote a $k$-basis for $A$.  Then we form the product $\frak{B} = \frak{A}^n$, note that it forms a $k$-basis for $B$, and let $G$ act on $\frak{B}$ via the inclusion $G \hookrightarrow \Sigma _n$.  For every orbit $\frak{o} = \{ \beta _1, \ldots , \beta _s \}$,  we let $t(\frak{o})$ denote the element $\beta_1 +  \cdots +\beta _s \in B^{G}$.  It is clear that the elements $t(\frak{o})$, as $\frak{o}$ ranges over the orbits of the $G$ action on $\frak{B}$, form a $k$-basis for $B^G$, and therefore that it suffices to prove that each element $t(\frak{o})$ admits a $p$-th root in $B^G$.  We suppose that we are given an orbit $\frak{o}$, with orbit representative $\beta \in \frak{o}$.  Let $G_\beta$ denote the stabilizer of $\beta$.  The group $G_{\beta}$ now acts on the set $\{ 1, \ldots , n \}$ via its inclusion into $\Sigma _n$, and we let 
$$ \{ 1 , \ldots , n \} =  S_1 \cup S_2 \cup \cdots  \cup S_l
$$
denote the orbit decomposition of this action.  Since $G_{\beta}$ stabilizes $\beta$, it is clear that there are elements $\{ \alpha _j \}_{j = 1}^l$, with $\alpha \in \frak{A}$, so that 
$$ \beta = \alpha _{j(1) }\otimes \alpha _{j(2)} \otimes \cdots \otimes \alpha _{j(n)}
$$
where $j(i)$ denotes the integer so that $i \in S_{j(i)}$.  For each $j = 1, \ldots, l$, we select an element $\overline{\alpha}_j$ so that  $(\overline{\alpha}_j)^p = \alpha _j$.  Also, select left coset representatives $\{ g_1, \ldots , g_k \}$ for the collection of left cosets $G/G_{\beta}$. The element 
$$ \sum _{m=1}^k g_m( \overline{\alpha}_{j(1)} \otimes \overline{\alpha}_{j(2)} \otimes \cdots \otimes \overline{\alpha}_{j(l)})
$$
is now clearly an element in $B^G$, and it is equally clearly a $p$-th root of $t(\frak{o})$, which gives us the claim.  We now proceed to the proof of the theorem.

 We consider first the case where $G \lhd G^{\prime}$ and where the index of $G$ in $G^{\prime}$ is $p$.   Consider any element $b \in I_{B^G}$.  We let $\overline{b}$ denote any $p$-th root of $b$ in $B^G$.  It is immediate  that $\overline{b} \in I_{B^G}$ as well.  Now construct the  polynomial 
$$ \varphi (X) = \prod _{\gamma  \in G^{\prime}/G} (X - \overline{b}^{\gamma})
$$
Note that the notation $\overline{b}^{\gamma}$ is well defined for $\gamma \in G^{\prime}/G$ because $\overline{b} \in B^G$.  
From the definition of $\varphi$, it is clear that $\varphi (\overline{b}) = 0$, and that $\varphi$  is of the form 
$$\varphi (X) = X^p + \sum _{i=1}^p (-1)^i c_i X^{p-i}
$$ 
where $c_i$ is the $i$-th elementary symmetric function in the elements $\overline{b}^{\gamma}$, as $\gamma $ varies over all the elements of $G^{\prime}/G$.  Note that $c_i \in  I_{B^{G^{\prime} }}$. 
Since $\overline{b}^p = b$, it follows that 
$$ b =   \sum _{i=1}^p (-1)^i c_i \overline{b}^{p-i}
$$
which gives the result in this case.  For the general case, we may form a sequence of subgroups 
$$  G  = G_0 \subseteq G_1 \subseteq G_2 \subseteq \cdots \subseteq G_{t-1} \subseteq G_t = G^{\prime}
$$
where each inclusion is the inclusion of a normal subgroup of index $p$.  This is an immediate consequence of the fact that $G^{\prime}$ is a $p$-group.  The result for the general case can now be obtained from the case of a normal subgroup via a downward induction on $t$.  We wish to prove that $I_{B^{G_0}} B^{G_t} = I_{B^{G_t}}$, and we suppose that we have proved that $I_{B^{G_{t-s}}}B^{G_t} =  I_{B^{G_t}}$.  We note that the case of a normal subgroup shows us that $I_{B^{G_{t-s-1}}} B^{G_{t-s}} = I_{B^{G_{t-s}}}$.  We  now obtain 
$$ I_{B^{G_t}} =  I_{B^{G_{t-s}}} B^{G_t}  = I_{B^{G_{t-s-1}}} B^{G_{t-s}} B^{G_t} =  I_{B^{G_{t-s-1}}} B^{G_t}
$$
The induction now proceeds until $s=t$ gives us the desired result.  
\end{Proof} 

We now  obtain the result which allows us to apply Corollary \ref{tootechnical}. 

\begin{Theorem}  \label{maindiagram} Let $k$ be a field of characteristic $p$, and $A$ an augmented  $k$-algebra in which every element has a $p$-th root.  Let $ G \subseteq G^{\prime} \subseteq \Sigma _n$ be an inclusion of $p$-groups, and let $ \displaystyle B = \mathop{\otimes} ^{n}_{k} A$ be equipped with the natural action of $\Sigma_n$ and therefore of $G$ and $G^{\prime}$. Also equip $B$ with the tensor product augmentation.  Suppose we have a diagram of augmented (to $k$) rings
$$
\begin{diagram}
\node{R^{\prime}} \arrow{s,t}{f} \arrow{e,t} {\pi ^{\prime}} \node{B^{G^{\prime}}} \arrow{s} \\
\node{R} \arrow{e,t}{\pi} \node{B^{G}}
\end{diagram} 
$$
in which $\pi $ and $\pi ^{\prime}$ are surjective, and so that the kernel $\frak{I}_R$ of $\pi$ satisfies $(\frak{I}_R)^N = \{ 0 \}$ for some positive integer $N$.  Then there exists an integer $M$ so that 
$(I_R)^M \subseteq  I_{R^{\prime}}\cdot  R$, where $I_R$ and $I_{R^{\prime}}$ denote the augmentation ideals of $R$ and $R^{\prime}$ respectively.

\end{Theorem} 
\begin{Proof}  We begin with $r \in I_R$, and apply $\pi$.  $\pi (r)$ is an element of $I_{B^G}$, and can therefore by Theorem \ref{general} be written as 
$$ \pi (r) = \sum _q i_q \beta _q
$$
where $i_q \in I_{B^{G^{\prime}}}$ and $\beta _q \in B^G$.   Because of the surjectivity of $\pi$ and $\pi ^{\prime}$, we can obtain elements  $\overline{i}_q \in R^{\prime} $ and $\overline{\beta} _q \in R$ so that $\pi (\overline{i}_q) = i_q $ and $\pi (\overline{\beta}_q) = \beta _q$.  Note that $\overline{i}_q \in I_{R^{\prime}}$.   Consider the element 
$$ \theta = \sum _{q} f(\overline{i}_q) \overline{\beta}_q  - r \in I_{R} 	$$
By construction  $\theta$ lies in  $\frak{I}_R$.  Consequently, by the hypothesis on $\frak{I}_R$, we know that $\theta ^N = 0$.  Let $M$ be any power of $p$ which is greater than $N$.  Then we have 
$$ 0 = \theta ^M = ( \sum _q f(\overline{i}_q) {\overline{\beta}_q})^M - r^M
$$
so
$$  r^M =  ( \sum _q f(\overline{i}_q) {\overline{\beta}_q})^M
$$
The right hand term is in $I_{R^{\prime}} \cdot R$, and since $M$ was chosen independently of $r$, we conclude that $(I_R) ^M \subseteq I_{R^{\prime}} \cdot R$, which was to be proved. 
\end{Proof}

We will use this result to study certain representation rings.  A first step is to study the behavior of the mod-$p$ representation rings of certain wreath products.  Recall that for any groups $G$ and $H$, with $G$ finite, the wreath product $G \wr H$ is the semidirect product $G \ltimes H^{\# (G)}$, where the $G$ action on $H^{\# (G)}$  is by permutation of coordinates.  For any group $G$, we will write $\frak{R}_p [G]$ for $R[G] \column{\otimes}{\Bbb{Z}} \Bbb{F}_p$.  

\begin{Proposition} \label{cycliccase}  Let $C_p$ denote the cyclic group of order $p$, where $p$ is a prime.  Then for any profinite group $K$, we construct $C_p \wr K$, with the inclusion $i: K^p \hookrightarrow C_p \wr K$.  The then following results hold.  
\begin{enumerate}
\item{  The restriction map $i^* : \frak{R}_p [C_p \wr K] \rightarrow \frak{R}_p [ K^p] $ is surjective onto the invariant subring $ \frak{R}_p [ K^p] ^{C_p}$.  }
\item{The kernel $\frak{I}$ of $i^*$ satisfies $\frak{I}^p = 0$. }
\end{enumerate} 
\end{Proposition} 
\begin{Proof}  Consider first the case where $K$ is finite.  Then the results of section 4.3 of \cite{james} show that an additive basis of $\frak{R}_p [C_p \wr K] $ can be constructed as follows.  Let $\mbox{Irr}(K)$ denote the set of irreducible representations of $K$.  Then we have the permutation action of $C_p$ on $\mbox{Irr}(K)^p$, which can then be decomposed into orbits, which either consist of one element or of $p$ elements.  Letting the orbit set be denoted by $\frak{O}$, we write $\frak{O} = \frak{O}_{(1)} \coprod \frak{O}_{(p)}$ for this decomposition.  A basis for $\frak{R}_p [C_p \wr K] $ is now in bijective correspondence with the set $\frak{O}_{(p)} \coprod \frak{O}_{(1)} \times \mbox{Irr}(C_p)$.  The correspondence is given by assigning to elements $\{ \rho _1, \ldots , \rho _p \}$  of $\frak{O}_{(p)}$ the induced representation $i_{K^p}^{C_p \wr K}(\rho _1)$, and on $ \frak{O}_{(p)} \times \mbox{Irr}(C_p)$ by assigning to elements $(\{\rho \} , \chi )$ the representation $\langle \rho, \chi \rangle$ whose underlying $K^p$ representation is $\otimes ^p \rho$, with the action by $C_p$ given by multiplication by the one-dimensional character $\chi$.  The restriction map $i^*$ is now clearly surjective onto the invariant subring $\frak{R}_p [ K^p] $, which is the first result.  The kernel ideal is  clearly spanned by differences of the form $\langle \rho, \chi \rangle - \langle \rho, \chi _0 \rangle $, where $\chi _0$ denotes the trivial character of $C_p$.  Such elements lie in
$I_{\frak{R}_p[C_p]}\cdot \frak{R}_p[C_p \wr K]$, where  $\frak{R}_p[C_p \wr K]$ is an $\frak{R}_p[C_p]$-algebra via the homomorphism induced by the group homomorphism  $C_p \wr K \rightarrow C_p$.  Since the augmenation ideal in $\frak{R}_p [C_p ]$ has trivial $p$-th power, the second result follows.  The extension to profinite groups is a straightforward direct limit argument.  
\end{Proof} 

Consider a $p$-Sylow subgroup $\Bbb{S}(p,n)$ of $\Sigma _{p^n}$.  It is well-known that this group can be described as an iterated wreath product 
$$ \Bbb{S}(p,n) = \wr^{i} C_p
$$
where $C_p$ denotes a cyclic group of order $p$.  We also  consider the semidirect product 
$\frak{S}(p,n) = \Bbb{S}(p,n)\ltimes \Bbb{Z}_p^{p^n}$.  

\begin{Proposition} \label{sylow}The image of the homomorphism $i^* : \frak{R}_p[ \frak{S}(p,n)] \rightarrow \frak{R}_p[\Bbb{Z}_p^{p^n}]$ is  the invariant subring $\frak{R}_p[\Bbb{Z}_p^{p^n}]^{\Bbb{S}(p,n)}$.  If we let $\frak{I}$ denote the kernel of the homomorphism, then $\frak{I}^{p^n} = \{ 0 \}$.  
\end{Proposition} 
\begin{Proof} Since $\frak{S}(p,n) \cong C_p \wr \frak{S}(p,n-1)$, a simple induction using Proposition \ref{cycliccase} gives the result.  
\end{Proof}

\begin{Corollary} \label{symmcomplete} 
The ring homomorphism  $\frak{R}_p[\Sigma _n \ltimes \Bbb{Z}_p^n ] \rightarrow \frak{R}_p[ \Bbb{Z}_p^n ]$ has image the invariant subring $ \frak{R}_p[ \Bbb{Z}_p^n ]^{\Sigma _n}$.  Further, its kernel $\frak{I}$ satisfies $\frak{I}^M = \{ 0 \}$, where $M$ is the smallest power of $p$ which is greater than $n$.  
\end{Corollary}
\begin{Proof}  We first observe that a $p$-Sylow subgroup $\Sigma _n(p)$ of $\Sigma _n$ is a product of $p$-Sylow subgroups of various subgroups $\Sigma _{p^i}$ of $\Sigma _n$, with the number of factors of the form $\Sigma _{p^i}$ being the coefficient of $p^i$ in the $p$-adic expansion of $n$.  Further, the actions of these subgroups are on disjoint blocks in a decomposition of the set $\{ 1, \ldots , n \}$ into blocks of size powers of $p$.  This means that the semidirect product  $\Sigma _n(p) \ltimes \Bbb{Z}_p^n$ decomposes into a product 
$$ \prod _{\alpha} \Sigma _{p^{i_{\alpha}}}(p^{i_{\alpha}}) \ltimes \Bbb{Z}_p^{p^{i_{\alpha}}}
$$ 
and therefore that $\frak{R}_p[ \Sigma _n(p) \ltimes \Bbb{Z}_p^n] $ decomposes into a tensor product 
$$ \bigotimes _{\alpha}\frak{R}_p [ \Sigma _{p^{i_{\alpha}}}(p^{i_{\alpha}}) \ltimes \Bbb{Z}_p^{p^{i_{\alpha}}}]
$$
Proposition \ref{sylow} applies to each of these factors, which means that the image of the ring homomorphism
$\frak{R}_p[ \Sigma _n(p) \ltimes \Bbb{Z}_p^n]  \rightarrow \frak{R}_p[ \Bbb{Z}_p^n] $ is the ring of invariants $\frak{R}_p[ \Bbb{Z}_p^n] ^{ \Sigma _n(p)}$.  Further, it is now immediate that the $M$-th power of the kernel of this homomorphism is trivial.  The result for $\frak{R}_p[\Sigma _n \ltimes \Bbb{Z}_p^n ] $ itself now follows by a simple transfer argument since the index of the $p$-Sylow subgroup in $\Sigma _n$ is prime to $p$.  
\end{Proof} 

\begin{Corollary} \label{fintech} Let $n$ be a positive integer, and let $k$ and $l$ be positive integers so that $k+l = n$.    The the homomorphism $\frak{R}_p[(\Sigma _k \times \Sigma_l )\ltimes \Bbb{Z}_p^n] \rightarrow \frak{R}_p [ \Bbb{Z}_p ^n]$ has image equal to the invariant subring $ \frak{R}_p [ \Bbb{Z}_p ^n]^{\Sigma _k \times \Sigma _l}$, and there is an integer $M$ so that the $M$-th power of the kernel ideal is $= \{ 0 \}$.  
\end{Corollary} 
\begin{Proof}  Immediate from Corollary \ref{symmcomplete}, since this example is simply a tensor product of two of the examples to which it applies.  
\end{Proof} 

\begin{Corollary} \label{reptech} Let $n, k$, and $l$ be as in the preceding corollary.  Then $\frak{R}_p[(\Sigma _k \times \Sigma_l )\ltimes \Bbb{Z}_p^n] $ becomes an  $\frak{R}_p[\Sigma _n \ltimes \Bbb{Z}_p^n ] $-algebra under the ring homomorphism induced by the inclusion $(\Sigma _k \times \Sigma_l )\ltimes \Bbb{Z}_p^n \rightarrow \Sigma _n \ltimes \Bbb{Z}_p^n$.  We let $I[-] \hookrightarrow \frak{R}_p[-]$ denote the augmentation ideal.  Then there is an integer $M$  so that 

$$ I [(\Sigma _k \times \Sigma_l )\ltimes \Bbb{Z}_p^n]^M \subseteq I[\Sigma _n \ltimes \Bbb{Z}_p^n ] \cdot \frak{R}_p[(\Sigma _k \times \Sigma_l )\ltimes \Bbb{Z}_p^n] 
$$
\end{Corollary} 
\begin{Proof}  This is a direct application of Theorem \ref{maindiagram}, with $R = \frak{R}_p[(\Sigma _k \times \Sigma_l )\ltimes \Bbb{Z}_p^n] $, $R^{\prime} = \frak{R}_p[\Sigma _n \ltimes \Bbb{Z}_p^n ] $, $A = \frak{R}_p[\Bbb{Z}_p]$, $f$ is the homomorphism induced by the inclusion 
$(\Sigma _k \times \Sigma_l )\ltimes \Bbb{Z}_p^n \hookrightarrow \Sigma _n \ltimes \Bbb{Z}_p^n $, and where $\pi $ and $\pi ^{\prime} $ are the homomorphisms induced by the inclusions $\Bbb{Z}_p ^n \hookrightarrow (\Sigma _k \times \Sigma_l )\ltimes \Bbb{Z}_p^n$ and $\Bbb{Z}_p ^n \hookrightarrow \Sigma _n\ltimes \Bbb{Z}_p^n$ respectively.  That the hypothesis of Theorem \ref{maindiagram} is satisfied is guaranteed by  Corollary \ref{fintech}. 
\end{Proof} 

\section{Equivariant $K$-theory for certain torus actions} \label{equivariant}

Let $\Bbb{A}^n(k) = Spec(k[x_1, x_2, \ldots , x_n ])$ be the $n$-dimensional affine space over $k$, with explicit choice of coordinates $x_i$.  Of course, Theorem 7 of \cite{quillen} tells us that $K(\Bbb{A}^n) \cong K(k)$, but it will turn out that a certain filtration of $K(\Bbb{A}^n) $ will be useful.  Let $S \subseteq \{ 1,2,\ldots ,n \}$ be any subset, and let ${\wp}(S) $ denote the prime ideal $ \Sigma_{s \in S} (x_s)$, and let ${\Pi}(S) $ denote the collection of prime ideals $\{ {\wp }(T)| T \supseteq S \}$.  We will  let $\Bbb{A}^n_S(k)$ denote the subscheme  $\Bbb{A}^n(k) - \bigcup _{s \in S} V(x_s)$, where 
 $V(x_s)$ denotes the hyperplane defined by $x_s$, and $\overline{\Bbb{A}}^n_S(k)$ denote the subscheme $\{ (x_1, \ldots ,x_n) | x_i =0 \mbox{ for } i \notin S \}$.  We also let $\mbox{Nil}(\Bbb{A}^n_S (k))$ denote the category of coherent $\Bbb{A}^n_S(k)$-modules on which the generators $\{ x_i \}_{i \notin  S}$ act nilpotently.   We let $\mbox{Mod} (n) = \mbox{Mod}(\Bbb{A}^n(k))$ denote the category of finitely generated $\Bbb{A}^n(k)$-modules.   For any $0 \leq i \leq n$, we let $\mbox{Mod}_i(n)$ denote the full subcategory of $\Bbb{A}^n$-modules $M$ for which the set of associated primes $Ass(M)$ is contained in the set

$${\Pi}_i = \bigcup _{\# (S) = i} {\Pi }(S)$$

See \cite{eisenbud} for a discussion of associated primes.  Similarly, for any set $S$ of cardinality $i$,  we define $\mbox{Mod}_{i+1}(n) \subseteq \mbox{Mod}_i(n;S)\subseteq \mbox{Mod}_i(n)$  to be the full subcategory on the  modules $M$ for which $Ass(M) \subseteq  \wp (S) \cup {\Pi}_{i+1}$.  We clearly have an increasing sequence of subcategories

$$\mbox{Mod}_n(n) \subseteq \mbox{Mod}_{n-1}(n) \subseteq \cdots \subseteq \mbox{Mod}_1(n) \subseteq \mbox{Mod}_0(n) = \mbox{Mod}(n)
$$

Further, each inclusion is the inclusion of an abelian subcategory, and it is clear from the discussion of quotient abelian categories in \cite{swan}  that the quotient $\mbox{Mod}_i (n)/ \mbox{Mod}_{i+1}(n) $ can be analyzed as follows.  For each $S$ with $\# (S) = i$, we may consider the inclusion $i_S : \mbox{Mod}_i(n;S) \hookrightarrow  \mbox{Mod}_i(n)$, and therefore the inclusion of subquotients $j_S: \mbox{Mod}_i(n;S) / \mbox{Mod}_{i+1}(n) \hookrightarrow  \mbox{Mod}_i (n)/ \mbox{Mod}_{i+1}(n) $.  We may then sum over all $S$ of cardinality $i$  and apply $K$-theory to obtain 

\begin{equation} \label{special}\bigvee_{\#(S) = i} j_S: \bigvee_{\#(S) = i} K(\mbox{Mod}_i(n;S)  / \mbox{Mod}_{i+1}(n))  \longrightarrow  K(\mbox{Mod}_i(n) / \mbox{Mod}_{i+1}(n))
\end{equation}

This functor is an equivalence since as in the proof of Theorem 5.4 of \cite{quillen} it can be checked that the underlying functor is an equivalence of categories.  Further, there is a natural equivalence of abelian categories 
$$ \mbox{Mod}_i(n;S)  / \mbox{Mod}_{i+1}(n) \simeq \mbox{Nil}(\Bbb{A}^n_{\overline{S}}(k))
$$

where $\overline{S} = \{ 1, \ldots , n \} -S$,  and hence  equivalences of $K$-theory spectra of these categories 
$$ K(\mbox{Mod}_i(n;S)  / \mbox{Mod}_{i+1}(n)) \simeq K(\mbox{Nil}(\Bbb{A}^n_{\overline{S}}(k)))
 \simeq K(\mbox{Mod}(\overline{\Bbb{A}}^n_{\overline{S}}(k)))
$$
The second of the equivalences  is proved as in \cite{quillen} by the devissage principle.







It is further clear from Section \ref{ktheory}  that all the spectra $$K(\mbox{Mod}_{i}(n)),K(\mbox{Mod}_i(n;S)), \mbox{ and }K(\mbox{Nil}(\Bbb{A}^n_S(k)))$$ are module spectra over the commutative  ring spectrum $K(\Bbb{A}^n(k))$,  that all the inclusion and restriction functors induce maps of $K(\Bbb{A}^n(k))$-module spectra, and therefore that the quotient spectra $K (\mbox{Mod}_{i}(n))/K(\mbox{Mod}_{i+1}(n))$ are also $K(\Bbb{A}^n(k))$-module spectra.  It is also clear that $K(\mbox{Nil}(\Bbb{A}^n_S(k)))$ is naturally a $K(\Bbb{A}^n_S(k))$-module spectrum, in a way which extends the $K(\Bbb{A}^n(k))$-module structure via the map of  ring spectra $K(\Bbb{A}^n(k)) \rightarrow K(\Bbb{A}^n_S(k))$.  Moreover, it is clear that the decomposition of spectra

$$ K(\mbox{Mod}_i(n)) / K(\mbox{Mod}_{i+1}(n)) \cong \bigvee_{\# (S) = i }  K(\mbox{Nil}(\Bbb{A}^n_{\overline{S}}(k)))
\cong \bigvee_{\# (S) = i }  K(\mbox{Mod}(\overline{\Bbb{A}}^n_{\overline{S}}(k)))
$$

given  by the homotopy inverse to the map $\bigvee_{\# (S) = i } K(j _s)$ from \ref{special} above  is also an equivalence of $K(\Bbb{A}^n(k))$-module spectra, where the $K(\Bbb{A}^n(k))$ action on the summand $K(\Bbb{A}^n_S(k))$ is the one obtained by restriction of scalars along the map of  ring spectra $K(\Bbb{A}^n(k)) \rightarrow K(\Bbb{A}^n_S(k))$.

We also wish to examine the behavior of the filtration and the decomposition of the subquotients under certain maps between affine spaces.  We let $X$ and $Y$ denote copies of the scheme $\Bbb{A}^n(k)$, and let $f^q : X \rightarrow Y$ be the power map given on coordinates by $f^q(x_1, \ldots , x_n) = (x_1^q, \ldots , x_n^q )$.  Let $X_S$ and $Y_S $ be the schemes obtained by removing the hyperplanes corresponding to the elements of $S$.  Let $\overline{X}_S \subseteq X_S $ and $\overline{Y}_S \subseteq Y_S $ denote the subschemes in which the variables $\{ x_i \}_{i \notin S}$ are equal  to zero.  Given any other Noetherian scheme $W$, we will also write $\mbox{Nil}(X_S;W)$ and $\mbox{Nil}(Y_S;W)$ for the full subcategory of modules on which the coordinates $\{ x_i \}_{i \notin S} $ act nilpotently.  We will write $\mbox{Mod}_i(X)$ and $\mbox{Mod}_i(Y)$ for the subcategories above in each of the two cases.  It is clear that these subcategories are respected by the restriction map along the map $f^q$.  We can now summarize the above discussion in the following Proposition. We also want to analyze the product of the affine spaces in question with a fixed scheme  $W$, and we will write $K(X;W)$ for $K(X \times W)$.

\begin{Proposition}  Let $X$ and $Y$ be copies of the affine scheme $\Bbb{A}^n(k)$, and let $f^q: X \rightarrow Y$ be the power map defined above. Let $W$ denote a fixed scheme.   The $K(X;W)$(resp. $K(Y;W))$-module spectra $K(\mbox{\em Mod}_i(X;W))\stackrel{defn}{=} K(X;W)^i$ and $K(\mbox{\em Mod}_i(Y;W))=K(Y;W)^i$ yield  filtrations of $K(X;W)$(resp. $K(Y:W ))$-module spectra, with $K(X;W)^0 = K(X;W)$ and $K(Y;W)^0 = K(Y;W)$.  The map $K(Y;W) \stackrel{K(f^q ;W)} \longrightarrow K(X;W)$ respects the filtration in the sense that there are evident maps $K(Y;W)^i \rightarrow K(X;W)^i$, so that the diagrams involving inclusion maps and the maps $K(f^q;W)$ all commute, yielding induced maps of subquotient $K(Y;W)$-module spectra $$K(Y;W)^i/K(Y;W)^{i+1} \rightarrow  K(X;W)^i/K(X;W)^{i+1}$$  
where the $K(Y;W)$-module structure on $K(X;W)^i/K(X;W)^{i+1}$ is obtained by restricting scalars along the map of  ring spectra $K(Y;W) \rightarrow K(X;W)$. 
The subquotients are themselves decomposed as 
$$ K(X;W)^i/K(X;W)^{i+1}  \cong \bigvee_{\# (S) = i} K(\mbox{\em Nil} (X_{\overline{S}};W)) \cong  \bigvee_{\# (S) = i} K( \overline{X}_{\overline{S}};W) $$ and $$K(Y;W)^i/K(Y;W)^{i+1} \cong \bigvee_{\# (S) = i} K(\mbox{\em Nil}(Y_{\overline{S}};W))  \cong \bigvee_{\# (S) = i} K( \overline{Y}_{\overline{S}};W)
$$
where the decomposition is as $K(Y;W)$-module spectra. The summands corresponding to $S$ for $X$ (resp. $Y$) are  equipped with a fixed  extension of the $K(X;W)$(resp. $K(Y;W)$)-module structure to a $K(X_{\overline{S}};W)$(resp. $K(Y_{\overline{S}};W)$)-module structure.  The subquotient map induced by $f_q$ respects this decomposition of $K(Y;W)$-module spectra, and on the summands corresponding to $S$, it respects the extended module structures in the sense that the map $K(\mbox{\em Nil} (Y_{\overline{S}};W))\rightarrow K(\mbox{\em Nil} (X_{\overline{S}};W))$ is a map of $K(Y_{\overline{S}};W)$-module spectra when $K(\mbox{\em Nil} (X_{\overline{S}};W))$ is regarded as a $K(Y_{\overline{S}};W)$-module spectrum via the  ring spectrum map  $K(Y_{\overline{S}};W) \rightarrow K(X_{\overline{S}};W)$.  
\end{Proposition}

We next begin the analysis of how this filtration works equivariantly.    Given $ X = \Bbb{A}^n(k)$ as above, with coordinates $\{ x_1, \ldots , x_n \}$,  there is an inclusion $\Sigma_n \ltimes (k^*)^n \hookrightarrow Aut(X)$.  Let $\rho : G \rightarrow \Sigma_n \ltimes (k^*)^n$ be any homomorphism, and let $\overline{\rho}: G \rightarrow \Sigma _n$ denote the corresponding permutation representation.  The  map $\rho$ now gives an action of the group $G$ on the scheme $X$, and  therefore on the module category $\mbox{Mod}(X)$.  It is clear that the filtration of $\mbox{Mod}(n)$ by the subcategories $\mbox{Mod}_i(n)$ is invariant under the action defined by $\rho$, and we therefore easily obtain a filtration of $\mbox{Mod}^G(X)$ by Serre subcategories $\mbox{Mod}_i^G(X)$.  We will have several categories  of $X$-modules which are closed under the action of $G$ or of subgroups of $G$.  For each such category, we can then denote its equivariant version with the superscript $G$, and define  its equivariant $K$-theory as the $K$-theory of the equivariant version of the category.  For instance, we will write $K^G(\mbox{Mod}_i(X))$ for $K(\mbox{Mod}^G_i(X))$.   We now identify the subquotients $\mbox{Mod}_i^G(X)/\mbox{Mod}_{i+1}^G(X)$.  We define  $X_i$ to be the scheme

$$ \coprod_{\#(S) = i}X_S
$$

We equip this scheme with a $G$-action as follows.  Let ${\cal P}_{i}$ denote the set of subsets of $\{ 1, \ldots , n \}$ of cardinality $n-i$, and let 
$$ {\cal P}_i =  \frak{o}_1 \cup \cdots \cup \frak{o}_t
$$
be an orbit decomposition, with orbit representatives $S_l \in \frak{o}_l$.  Let $G_l$ denote the stabilizer of $S_l $ in $G$.  Then we may construct the scheme with $G$-action
$$  \coprod _l G \column{\times}{G_l} X_{S_l}
$$
and define (uniquely)  an equivariant map 

$$  \coprod _l G \column{\times}{G_l}X_{S_l} \rightarrow X
$$
by extensions of each of the $G_l$-equivariant inclusions $  X_{S_l} \rightarrow X$.  This map clearly gives an identification of $\coprod _l G \column{\times}{G_l} X_{S_l} $ with $X_{n-i}$, and we have therefore defined a $G$-action on $X_{n-i}$. 
Further, we  let $\overline{X}_{n-i} $ denote the subscheme $\coprod _{\# (S) = n-i} \overline{X}_S$, and define a subcategory $\mbox{Nil}(X_{n-i})$ of the category of all finitely generated $X_{n-i}$-modules as follows.  A finitely generated $X_{n-i}$-module $M$ is simply a choice of $X_S$-module $M_S$ for each $S \subseteq \{ 1, \ldots , n \}$ of cardinality $n-i$.  $\mbox{Nil}(X_{n-i})$ is defined to be the category of all $M$ so that $M_S \in \mbox{Nil}(X_S)$ for all such $S$.

By an argument again entirely parallel to that in the proof of Theorem 5.4 of \cite{quillen}, we have that 

$$ K(\mbox{Mod}_i^G(X)/\mbox{Mod}_{i +1}^G(X)) \cong K^G(\mbox{Nil}(X_{n-i}))
$$

For each orbit $\frak{o}_l$, we let $X_{\frak{o}_l}$ denote the corresponding subscheme of $X_i$, and similarly for $\overline{X}_{\frak{o}_l}$.  The category $\mbox{Nil}$ is  of course closed under the $G$-action, and there is an evident equivariant decomposition of categories

$$ \mbox{Nil}(X_{n-i}) \cong \prod_l\mbox{Nil}(X_ {\frak{o}_l})
$$
and therefore a decomposition of $K$-theory spectra 

$$ K(\mbox{Mod}^G(X)_i/\mbox{Mod}^G(X)_{i +1})  \cong 
\bigvee _{l=1}^t K^G(\mbox{Nil} (X_{\frak{o}_l})) \cong  \bigvee _{l=1}^t K^G( \overline{X} _{\frak{o}_l})
$$

Given a Noetherian scheme $W$ with $G$-action, one can by analogy with the non-equivariant case construct  categories $\mbox{Mod}_i^G(X,W)$ and $\mbox{Nil} ^G(X_i;W)$ as well as the subcategories  $\mbox{Nil}^G (X_{\frak{o}_l};W)$, along with the corresponding $K$-theory spectra, and obtain a parallel decomposition

$$ K(\mbox{Mod}^G_i(X;W)/\mbox{Mod}_{i +1}^G(X;W))  \cong 
\bigvee _{l=1}^t K^G(\mbox{Nil} (X_{\frak{o}_l};W)) \cong  \bigvee _{l=1}^t K^G( \overline{X} _{\frak{o}_l};W)
$$

We now wish to describe the behavior under the power maps described above.  We let $X$ and $Y$ denote copies of $\Bbb{A}^n(k)$ equipped with actions by a profinite  $p$-group $G$ via homomorphisms $\rho _X, \rho _Y : G \rightarrow \Sigma _n \ltimes (k^*)^n$ as above, and we assume that the power map $f^q : X \rightarrow Y$ is given as above, and further that the actions of $G$ on $X$ and $Y$ are compatible with this map in the sense that $f^q$ is $G$-equivariant when $X$ and $Y$ are equipped with the actions $\rho _X$ and $\rho _Y$, respectively.  It is immediate from this assumption that the composites $G \stackrel{\rho _X}{\rightarrow} \Sigma _n \ltimes (k^*)^n \rightarrow \Sigma _n$ and $G \stackrel{\rho _Y }{\rightarrow} \Sigma _n \ltimes (k^*)^n \rightarrow \Sigma _n$ are identical.   We will let $\mbox{Mod}_i^G(Y;W )$, $\mbox{Nil}^G(Y_{\frak{o}_l};W)$, $Y_{\frak{o}_l}$, and $\overline{Y}_{\frak{o}_l}$ be defined by analogy with $\mbox{Mod}_i^G(X;W )$, $\mbox{Nil}^G(X_{\frak{o}_l};W)$, $X_{\frak{o}_l}$, and $\overline{X}_{\frak{o}_l}$

\begin{Proposition} \label{indsystem}  Let $X$ and $Y$ and  $f^q: X \rightarrow Y$ as above, and let $W$ denote any fixed scheme with $G$-action.  The $K^G(X;W)$(resp. $K^G(Y;W))$-module spectra $K^G(\mbox{\em Mod}_i(X;W))\stackrel{defn}{=} K^G(X;W)^i$ and $K^G(\mbox{\em Mod}_i(Y;W))=K^G(Y;W)^i$ yield  filtrations of $K^G(X;W)$(resp. $K^G(Y;W))$-module spectra, with $K^G(X;W)^0 = K^G(X;W)$ and $K^G(Y;W)^0 = K^G(Y;W)$.  The map $K^G(Y;W) \stackrel{K^G(f^q;W)} \longrightarrow K^G(X;W)$ respects the filtration in the sense that there are evident maps $K^G(Y;W)^i \rightarrow K^G(X;W)^i$, so that the diagrams involving inclusion maps and the maps $K^G(f^q;W)$ all commute, yielding induced maps of subquotient $K^G(Y;W)$-module spectra $$K^G(Y;W)^i/K^G(Y;W)^{i+1} \rightarrow  K^G(X;W)^i/K^G(X;W)^{i+1}$$  
where the $K^G(Y;W)$-module structure on $K^G(X;W)^i/K^G(X;W)^{i+1}$ is obtained by restricting scalars along the map of  ring spectra $K^G(Y;W) \rightarrow K^G(X;W)$. 
The subquotients are themselves decomposed as

$$ K^G(X;W)^i/K^G(X;W)^{i+1} \cong \bigvee_{l=1}^t K^G(\mbox{\em Nil}(X_{\frak{o}_l};W))
\cong \bigvee_{l=1}^t  K^{G_l}(\mbox{\em Nil}(X_{S_l};W))$$

and

$$K^G(Y;W)^i/K^G(Y;W)^{i+1} \cong \bigvee_{l=1}^t K^G(\mbox{\em Nil}(Y_{\frak{o}_l };W))
\cong \bigvee_{l=1} ^tK^{G_l}(\mbox{\em Nil}(Y_{S_l};W))
$$

where the decomposition is as $K^G(Y;W)$-module spectra. The summands corresponding to $\frak{o}_l$ for $X$ (resp. $Y$) are  equipped with a fixed  extension of the $K^G(X;W)$(resp. $K^G(Y;W)$)-module structure to a $K^G(X_{\frak{o}_l};W)$(resp. $K^G(Y_{\frak{o}_l};W)$)-module structure.  The subquotient map induced by $f_q$ respects this decomposition of $K^G(Y;W)$-module spectra, and on the summand corresponding to ${\frak{o}}$, it respects the extended module structures in the sense that the map $K^G(\mbox{\em Nil}(Y_{\frak{o}_l};W)) \rightarrow K^G(\mbox{\em Nil}(X_{\frak{o}_l};W))$ is a map of $K^G(Y_{\frak{o}};W)$-module spectra when $K^G(\mbox{\em Nil}(X_{\frak{o}_l};W))$ is regarded as a $K^G(Y_{\frak{o}_l};W)$-module spectrum via the  ring spectrum map  $K^G(Y_{\frak{o}_l};W) \rightarrow K^G(X_{\frak{o}_l};W)$. Finally each of the module spectra $K^G(\mbox{\em Nil}(Y_{\frak{o}_l};W))$ (respectively  $K^G(\mbox{\em Nil}(X_{\frak{o}_l});W)$ ) are free and cyclic  as $K^G(Y_{\frak{o}_l};W)$- (respectively $K^G(Y_{\frak{o}_l};W)$)-module spectra, and the map $K^{G_{\frak{o}_l}}(Y_{S_{\frak{o}_l}};W) \rightarrow K^{G_{\frak{o}_l}}(X_{S_{\frak{o}_l}};W)$   induced by $f^q$  is clearly equal to the map induced by the power map $ X_{S_{\frak{o}_l} }\rightarrow Y_{S_{\frak{o}_l}}$.  
\end{Proposition} 
\begin{Proof}
Everything in the proposition is clear from the previous discussion perhaps the free and cyclic properties of the modules $K^G(\mbox{Nil}(Y_{\frak{o}_l};W))$ and $K^G(\mbox{Nil}(X_{\frak{o}_l};W))$.  For $X$, this follows from the following commutative diagram.  
$$
\begin{diagram}
\node{K^G(\overline{X}_{\frak{o}_l};W)} \arrow{se} \arrow{s} \node{} \\
\node{K^G(X_{\frak{o}_l},W)} \arrow{e} \node{K^G(\mbox{Nil} (X_{\frak{o}_l};W))}
\end{diagram}
$$
The left hand vertical map is the map of  ring spectra induced from the evident projection of $G$-schemes $X_{\frak{o}_l} \rightarrow \overline{X}_{\frak{o}_l}$.  It induces an equivalence on $K$-theory spectra from the homotopy property for $K^G$-theory.  The diagonal arrow, on the other hand, is induced by the inclusion of the subcategory of $X_{\frak{o}_l}$-modules on which the nilpotence condition is strengthened to one which says that the subset of variables which are required to act nilpotently to be an object of $\mbox{Nil} (X_{\frak{o}_l};W)$ are instead required to act identically by zero.  This inclusion is also a weak equivalence due to the devissage theorem.  It follows that the bottom horizontal arrow is also an equivalence.  
\end{Proof}

Finally, we record how the decompositions behave under change of groups homomorphisms.  
Of course, when we have a group homomorphism $f: G \rightarrow H$ and a smooth Noetherian scheme $W$ with $H$-action, there is a naturally defined map of commutative  ring spectra $K^f(X;W) : K^H(X;W) \rightarrow K^G(X;W)$.  

\begin{Proposition} \label{changegroups} 
Let $f: G \rightarrow H$ be a group homomorphism, and suppose we have an action of $H$ on $ X = \Bbb{A}^n(k)$ via a homomorphism $\rho : H \rightarrow \Sigma _n \ltimes (k^*)$.  Then there are naturally defined maps of module spectra 
$K^f(X,W)^i : K^H(X,W)^i \rightarrow K^G(X,W)^i$, which commute with  all the relevant inclusion maps between them.  Consequently, there are induced maps between subquotients.  We then have the $H$ and $G$-actions on the sets ${\cal P}_i$ of subsets of $\{ 1, \ldots, n \}$ of cardinality $=i$, and orbit sets $\Omega _i^G$ and $\Omega _i^H$, with the map $\Omega _i^f: \Omega _i^G \rightarrow \Omega _i^H$.  For any $\frak{o} \in \Omega _i^H$, there are  natural maps

$$ K^H(X_{\frak{o}};W)\rightarrow \prod _{\omega \in (\Omega_i^f)^{-1} (\frak{o})}
K^G(X_{\omega};W)
$$
and 
$$ K^H(\mbox{\em Nil}(X_{\frak{o}};W)) \rightarrow \prod _{\omega \in (\Omega_i^f)^{-1} (\frak{o})}
K^G(\mbox{\em Nil}(X_{\omega};W))
$$
of module spectra induced by $f$.  Finally, the diagrams 
$$  
\begin{diagram}
\node{K^H(X_{\frak{o}};W)} \arrow{s,t}{\sim} \arrow{e} \node{\prod _{\omega \in (\Omega_i^f)^{-1} \frak{o}}
K^G(X_{\omega};W)} \arrow{s,b}{\sim} 
\\
\node{K^H(\mbox{\em Nil}(X_{\frak{o}};W))} \arrow{e} \node{\prod _{\omega \in (\Omega_i^f)^{-1} \frak{o}}
K^G(\mbox{\em Nil}(X_{\omega};W))}
\end{diagram} 
$$
commute.  
\end{Proposition} 

\section{Equivariant $K$-theory of ${\frak{G}}$-representations}

We  now carry out the analysis of the equivariant $K$-theory of monomial ${\frak{G}}$-representations.  We assume that the characteristic of $k$ is not equal to $p$.  A  ${\frak{G}}$-representation $\rho$ of a profinite $p$-group $ G$ is an action  on the special pro-group scheme

$$ {\cal G}= \mbox{ \hspace{.1cm}} (\cdots \rightarrow   \Bbb{G}_m(k)^{n,\rho} \stackrel{f^p}{\rightarrow}  \Bbb{G}_m(k)^{n,\rho} \stackrel{f^p}{\rightarrow}  \Bbb{G}_m(k)^{n,\rho}    \stackrel{f^p}{\rightarrow}  \Bbb{G}_m(k)^{n,\rho}  )
$$ 

for some $n$, where the superscript $\rho$ is a label indicating the group action, where $f^p$ is the power map defined in Section \ref{equivariant}.   We recall that the representation was constructed as  a homomorphism to the subgroup of automorphisms of the tower of schemes above given by $GL_n(\Bbb{Z}) \ltimes \Bbb{Z}_p^n \subseteq GL_n(\Bbb{Z}) \ltimes (
\column{lim}{\leftarrow}k^*)^n$.  
Because the action is monomial, it is clear that the action extends to the special pro-scheme 

$${\cal A}= \mbox{ \hspace{.1cm}}(\cdots \rightarrow \Bbb{A}_n^{\rho}(k)  \stackrel{\times p}{\rightarrow}  \Bbb{A}_n^{\rho}(k) \stackrel{\times p }{\rightarrow}   \Bbb{A}_n^{\rho}(k)  \stackrel{\times p }{\rightarrow}  \Bbb{A}_n^{\rho}(k)   \stackrel{\times p }{\rightarrow}  \Bbb{A}_n^{\rho}(k) )
$$

where $\times p :  \Bbb{A}_n(k)  \rightarrow  \Bbb{A}_n(k) $ denotes the map which raises each coordinate to the $p$-th power, and again the superscript $\rho$ indicates the action.  We will write ${\cal A}^r$ and ${\cal G}^r$ for the $r$-th term in the corresponding pro-schemes.  Since the $G$-action  attached to $\rho$  is determined by a choice of group action for each individual $\Bbb{G}_m(k)^{n,\rho}$, together with compatibilities as described in Definition \ref{creps}, we obtain a special Ind-spectrum

$$  K^G(\Bbb{G}_m(k)^{n, \rho }) \rightarrow  K^G(\Bbb{G}_m(k)^{n, \rho }) \rightarrow  K^G(\Bbb{G}_m(k)^{n, \rho }) \rightarrow  K^G(\Bbb{G}_m(k)^{n, \rho }) \rightarrow \cdots 
$$

which we denote by $K^G({\cal G})$.   We define $K^G({\cal A})$ similarly.  We let $W$ denote any scheme with a continuous action of the profinite group $G$, and define $K^G({\cal A};W)$ and $K^G({\cal G};W)$ as in Section \ref{equivariant}.  It is clear that there is a natural map $\lambda: K^G({\cal A};W) \rightarrow K^G({\cal G};W)$.  It is also clear from the existence of the  maps ${\cal A} \rightarrow *$ and ${\cal G} \rightarrow *$ that $K^G({\cal A};W)$ and $K^G({\cal G};W)$ are both special Ind-module spectra over the commutative  ring spectrum $K^G(*)$.  The ring spectrum $K^G(*)$ is equipped with a natural augmentation of commutative  ring spectra $\varepsilon: K^G(*) \rightarrow \Bbb{H}(\Bbb{F}_p)$, where $ \Bbb{H}(\Bbb{F}_p)$ denotes the mod $p$ Eilenberg-MacLane spectrum, regarded as a commutative ring spectrum in the sense that any ordinary commutative ring  is (see \cite{completion}).  The goal of this section is to prove that the  map 

$$  \Bbb{H}_p \column{\wedge}{K^G(*)}(\column{hocolim}{\rightarrow} K^G({\cal A};W))\rightarrow    \Bbb{H}_p \column{\wedge}{K^G(*)}(\column{hocolim}{\rightarrow }K^G({\cal G};W))
$$

induced by $\lambda$ is a weak equivalence of spectra, where the smash products are derived smash products, i.e. where the homotopy colimit spectra are replaces by cofibrant objects in the module category over $K^G(*)$.  In order to do this, we must construct a filtration on $K^G({\cal A};W)$ by special Ind-$K^G(*)$-module spectra.  We recall the filtration of categories 

$$\mbox{Mod}_n (n)\subseteq \mbox{Mod}_{n-1}(n) \subseteq \cdots \subseteq \mbox{Mod}_1(n) \subseteq \mbox{Mod}_0(n) = \mbox{Mod}(n)
$$

from Section \ref{equivariant}.   It is clearly preserved under the power maps defining ${\cal A}$, and we may therefore speak of the special Ind-abelian category $\mbox{Mod}_i({\cal A};W)$.  It is also clear from the definitions that $K^G({\cal G};W) \cong K^G(\mbox{Mod}_n({\cal A};W))$ and that $K^G({\cal A};W) \cong K^G(\mbox{Mod}_0({\cal A};W))$.  Therefore, it suffices to prove that 
$$ \Bbb{H}_p \column{\wedge}{K^G(*)}(\column{hocolim}{\rightarrow} (K^G(\mbox{Mod}_i({\cal A};W)) /K^G(\mbox{Mod}_{i+1}({\cal A};W))) )\cong *$$
 for $0 \leq i \leq n-1$, since the operation $ \Bbb{H}_p \column{\wedge}{K^G(*)}( - )$ respects cofibration sequences of $K^G(*)$-modules, and since homotopy colimits of spectra also preserve cofibration sequences. 
 
 Let $ \sigma: G \rightarrow \Sigma _n$ denote the permutation representation associated to $\rho$.  As in Section \ref{equivariant}, we consider an orbit decomposition of the set ${\cal P}_i$ of subsets of $\{ 1, \ldots , n \} = \{ \frak{o}_1, \ldots , \frak{o}_t \}$.  From Proposition \ref{indsystem}, it is clear that we may construct a decomposition 
 
 $$K^G(\mbox{Mod}_i({\cal A};W)) /K^G(\mbox{Mod}_{i+1}({\cal A};W))
 \cong \bigvee _{l=1}^t K^G(\mbox{Nil}({\cal A}_{\frak{o}_l};W))
 $$
 
 where $\mbox{Nil}({\cal A}_{\frak{o}_l} ;W)$ is the special Ind-spectrum constructed from the  summands for the various categories $\mbox{Nil}({\cal A}^r_{\frak{o}_l},W)$.   In order to analyze these special Ind-spectra, we will some concepts. 
 
 \begin{Definition} Let $\Bbb{H}$ denote a (-1)-connected commutative  ring spectrum.  Let $Alg/\Bbb{H}$ denote the category of (-1) connected commutative ring spectra over $\Bbb{H}$, i.e. equipped with a homomorphism $A \rightarrow \Bbb{H}$ of commutative  ring spectra.  By an $s\frak{I}$-algebra over $\Bbb{H}$, we will mean a special Ind-object in the category $Alg/\Bbb{H}$.  We note that the category {\em sInd($\underline{\mbox{\em Spectra}}$)} admits a coherently associative and commutative smash product $\wedge$, by simply applying the smash product objectwise to the functor defining the special Ind object.  By a module over an $s\frak{I}$-algebra $(A, \varphi: A \rightarrow \Bbb{H})$, we will mean an sInd-spectrum $M$ together with a map $A \wedge M \rightarrow M$ satisfying the usual conditions defining a module.  Objectwise, this means that if we are given a special Ind-object $\{ R_i \}_{i\geq 0}$ in the category of  commutative ring spectra, then a module over it is a family of spectra $\{ M_i \}_{i \geq 0}$, together with maps $\theta _i : M_i \rightarrow M_{i+1}$ with each $M_i$ equipped with an $R_i$-module structure, and so that each map $\theta _i$ is a map of $R_i$-modules if $M_{i+1} $ is regarded as an $R_i$-module via restriction of scalars along the homomorphism $R_i \rightarrow R_{i+1}$ of commutative  ring spectra. A module over $\{ R_i \}$  is said to be {\em quasi-cyclic} if each $R_i$-module $M_i$ is weakly equivalent to a free module of rank one over $R_i$.  A quasi-cyclic module $\{ M_i \}$ over $\{ R_i \}$ is {\em based} if we are given specific equivalences of $R_i$-modules $R_i \rightarrow M_i$.  
  \end{Definition} 
  
  The following are readily verified. 
  
  \begin{Proposition} Let $\underline{{\cal R}} $ and $\underline{{\cal R}^{\prime}}$ be $s\frak{I}$-algebras, and let  $f: \underline{{\cal R}} \rightarrow \underline{{\cal R}^{\prime}}$ be a homomorphism of $s\frak{I}$-algebras, then there is an extension of scalars functor $\underline{{\cal R}^{\prime}} \column{\wedge}{\underline{{\cal R}}} -$, taking modules over $\underline{{\cal R}} $  to modules over $\underline{{\cal R}^{\prime}}$.   It has the usual universality property of this kind of construction, namely if we are given an $\underline{{\cal R}}$-module $\underline{{\cal M}} $, an $\underline{{\cal R}^{\prime}}$-module $\underline{{\cal M}^{\prime}} $, and a homomorphism $g:\underline{{\cal M}} \rightarrow \underline{{\cal M}^{\prime}}$ of 
  $\underline{{\cal R}}$-modules, then there is a unique extension of $g$ to a homomorphism of $\underline{{\cal R}^{\prime}}$-modules $\overline{g}: \underline{{\cal R}^{\prime}} \column{\wedge}{\underline{{\cal R}}}\underline{{\cal M}}\rightarrow \underline{{\cal M}^{\prime}}$ which restricts to $g$ on $M$. If $\underline{\cal {M}}$ is a quasi-cyclic $\underline{{\cal R}}$-module, then $\underline{{\cal R}^{\prime}} \column{\wedge}{\underline{{\cal R}}}\underline{{\cal M}}$ is a quasi-cyclic  $\underline{{\cal R}^{\prime}}$-module.\end{Proposition} 
  
  \begin{Proposition} \label{inflate} 
   Let $\underline{{\cal R}} $ and $\underline{{\cal R}^{\prime}}$ be $s\frak{I}$-algebras, and let  $f: \underline{{\cal R}} \rightarrow \underline{{\cal R}^{\prime}}$ be a homomorphism of $s\frak{I}$-algebras,  Let $\underline{{\cal M}} $ (respectively $\underline{{\cal M}^{\prime}} $) be a quasi-cyclic $\underline{{\cal R}} $ (respectively $\underline{{\cal R}^{\prime}}$)-module, and suppose we are given a homomorphism $h: \underline{{\cal M}} \rightarrow  \underline{{\cal M}^{\prime}} $ of $\underline{{\cal R}} $-modules.  Suppose further that we are given based structures for $\underline{{\cal M}} $ and $ \underline{{\cal M}^{\prime}} $, i.e. families of maps of spectra  $\{ \psi _i : R_i \rightarrow M_i \}$ and $\{ \psi ^{\prime}: R_i^{\prime} \rightarrow M_i^{\prime}\}$,  so that $\psi _i $ is an equivalence of $R_i$-modules, and $\psi ^{\prime}$ is an equivalence of $R_i^{\prime}$-modules.  If $h$ is compatible with the based structures in the sense that the diagrams
   $$
   \begin{diagram}
    \node{R_i} \arrow{s,t}{f_i}  \arrow{e,t}{\psi _i}  \node{M_i} \arrow{s,b}{h_i}  \\
   \node{R_i^{\prime}} \arrow{e,t}{\psi ^{\prime}} \node{M_i^{\prime}}
   \end{diagram}  
   $$
   commute, then the natural extension  of scalars  map $ \underline{{\cal R}^{\prime}} \column{\wedge}{\underline{{\cal R}} } \underline{{\cal M}} \rightarrow \underline{{\cal M}^{\prime}} $ is an equivalence of $\underline{{\cal R}^{\prime}}$-modules.  
  \end{Proposition} 
  \begin{Proof}  This is straightforward, and we omit it.  It is the natural extension of the corresponding fact for free cyclic modules over an ordinary ring.  
  \end{Proof} 

{\bf Example:}  The main example for us comes from the previous section.  The system
$$K^G({\cal A}_{\frak{o}_l};W) = \{  K^G({\cal A}^0_{\frak{o}_l};W) \rightarrow K^G({\cal A}^1_{\frak{o}_l};W) \rightarrow K^G({\cal A}^2_{\frak{o}_l};W) \rightarrow K^G({\cal A}^3_{\frak{o}_l};W) \rightarrow \cdots
\} $$
is an $s\frak{I}$-algebra with ground algebra  $K^G(*)$ augmented over $\Bbb{H}_p$, the mod-$p$ Eilenberg-MacLane spectrum.  
$$K^G(\mbox{Nil}({\cal A}_{\frak{o}_l};W)) = \{   K^G(\mbox{Nil}({\cal A}^0_{\frak{o}_l};W)) \rightarrow K^G(\mbox{Nil}({\cal A}^1_{\frak{o}_l};W)) \rightarrow K^G(\mbox{Nil}({\cal A}^2_{\frak{o}_l};W)) \rightarrow \cdots\}
$$
is a quasi-cyclic module over $K^G({\cal A}_{\frak{o}_l};W)$.  These statements follow directly from Proposition \ref{indsystem}.


\begin{Proposition} \label{msystemtech}
Let $A\rightarrow \Bbb{H}$ be a homomorphism of  $(-1)$-connected commutative ring spectra.   Let $I $ denote the kernel of the homomorphism $\pi _0 A \rightarrow\pi _0 \Bbb{H}$.  Let $\frak{B}$ and $\frak{C}$ denote $s\frak{I}$-algebras over $\Bbb{H}$. Let $f: \frak{B} \rightarrow \frak{C}$ be a homomorphism of $s{\cal I}$-algebras over $\Bbb{H}$.   Suppose further that we are given a based quasi-cyclic $\frak{B}$-module $\frak{M}$, where the based structure consists of equivalences of $B_i$-modules $\psi _i : B_i \rightarrow M_i$.   Then the based structure on $\frak{M}$ determines a family of elements ${\cal S}(\frak{M}) = \{ \alpha _i  \in \pi _0 B_i \}_{i \geq 1}$ via the equation $\pi _0(\theta _i \compcirc  \psi _i)(1)  = \alpha _{i+1}\pi _0 (\psi _{i+1})(1)$, where the maps $\theta _i: M_i \rightarrow M_{i+1}$ are the structure maps for $\frak{M}$.  
We now have the following results. 
\begin{enumerate}
\item{${\cal S}(\frak{C}\column{\wedge}{\frak{B}} \frak{M} ) = \pi _0(f) {\cal S}( \frak{M})$, where $\pi _0(f) {\cal S}( \frak{M}) $ denotes the set $\{ \pi _0 (f_i) ( \alpha _i) \}_{i \geq 1}$.  }
\item{Let the maps $\rho_i : B_i \rightarrow B_{i+1}$ denote the structure maps for $\frak{B}$, and let $r_i^j :B_i \rightarrow B_j$ denote the composite $\rho _{j-1}\rho_{j-2} \cdots \rho _i$.  Suppose that there is an integer $M$ so that all products 
$$\alpha _{i+M} \cdot r^{i+M}_{i+ M-1}(\alpha_{i+M-1})  \cdot r_{i+ M-2}^{i+M}(\alpha_{i+M-2}) \cdots r_{i}^{i+M}(\alpha_{i}) 
$$
lie in $I \cdot \pi _0 B_{i+M}$, for all values of $i$.  
Then $\Bbb{H} \column{\wedge}{A}( \column{hocolim}{\rightarrow }\frak{M}) \simeq *$.  }
\end{enumerate}
\end{Proposition}
\begin{Proof}  The first statement is immediate.  The second follows directly from Proposition \ref{localize}. 
\end{Proof}

\begin{Corollary} Let $\{ \alpha _i \}_{i \geq 1} $ be the set 
$$   {\cal S}( K^G(\mbox{Nil}({\cal A}^0_{\frak{o}_l}) ) \rightarrow 
K^G(\mbox{Nil}({\cal A}^1_{\frak{o}_l}) ) \rightarrow K^G(\mbox{Nil}({\cal A}^2_{\frak{o}_l}) ) \rightarrow \cdots )
$$
which makes sense since the system in question is a based quasi-cyclic module  over the $s\frak{I}$-algebra
$ K^G({\cal A}_{\frak{o}_l}) $ by Proposition \ref{indsystem}. 
Of course,  $\pi _0 K^G(*) \cong R[G]$, $\pi _0 \Bbb{H}_p \cong \Bbb{F}_p$, and the augmentation induces the mod $p$-reduction of the augmentation on $R[G]$, whose kernel we denote $I_p[G]$. As in Proposition \ref{msystemtech}, let $r_i^j: K^G(\mbox{Nil}({\cal A}^i_{\frak{o}_l}) ) \rightarrow K^G(\mbox{Nil}({\cal A}^j_{\frak{o}_l}) ) $ be the structure map for each $i$ and $j$.  If the collection $\{ \alpha _i \}$ has the property that there exists an integer $M$ so that  
$$\alpha _{i+M} \cdot r_{i+M-1}^{i+M}(\alpha_{i+M-1})  \cdot r_{i+ M-2}^{i+M} (\alpha_{i+M-2})  \cdots r_i^{i+M}(\alpha_{i})   \in I_p[G] \cdot \pi _0 K^G({\cal A}^{i+M}_{\frak{o}_l})$$ for all $i$, then $$\Bbb{H}_p \column{\wedge}{K^G(*)} (\column{hocolim}{\rightarrow} K^G(\mbox{\em Nil}({\cal A}_{\frak{o}_l},W) )  ) \simeq *$$ for any smooth Noetherian $G$-scheme $W$.
\end{Corollary} 
\begin{Proof}  For $i \geq 0$, we let $p_i : K^G({\cal A}^i_{\frak{o}_l}) \rightarrow K^G({\cal A}^i_{\frak{o}_l};W) $ be the map induced by the projection $W \rightarrow *$.   By statement (1) of  Proposition \ref{msystemtech}, we have that 
$${\cal S} (K^G(\mbox{Nil}({\cal A}^0_{\frak{o}_l},W) ) \rightarrow 
K^G(\mbox{Nil}({\cal A}^1_{\frak{o}_l},W) ) \rightarrow K^G(\mbox{Nil}({\cal A}^2_{\frak{o}_l},W) ) \rightarrow \cdots)
$$
$$  = \{ \pi _0(p_i) (\alpha _i) \}_{i \geq 1}  
$$
If 
$$\alpha _{i+M} \cdot r_{i+M-1}^{i+M}(\alpha_{i+M-1})  \cdot r_{i+M-2}^{i+M}(\alpha_{i+M-2}) \cdots r_{i}^{i+M} (\alpha_{i})   \in I_p[G] \cdot \pi _0 K^G({\cal A}^{i+M}_{\frak{o}_l})$$
it follows immediately that 
$$ \pi _0(p_{i+M}) (\alpha _{i+M})\cdot  r_{i+M-1}^{i+M}( \pi _0(p_{i+M-1})(\alpha_{i+M-1})) \cdots r_{i}^{i+M}( \pi _0(p_{i})(\alpha_{i})) $$
$$\in I_p[G] \cdot \pi _0  K^G({\cal A}^{i+M}_{\frak{o}_l};W)$$which gives the result by statement (2) of Proposition \ref{msystemtech}. 
\end{Proof}

We now analyze the directed system  $K^G(\mbox{Nil}({\cal A}_{\frak{o}_l}) ) $.
Recall that $\frak{o}_l$ is an orbit under the action of $G$ on the collection of non-empty subsets of $\{ 1, \ldots , n \}$.  Let $S_l \in \frak{o}_l$ be an orbit representative.   Let $\Sigma _{S_l} \subseteq \Sigma _n$ denote the stabilizer of $S_l$ in $\Sigma _n$.  It is a subgroup of the form $\Sigma _s \times \Sigma _t \subseteq \Sigma _n$, with $s+t = n$.  Then the actions of the stabilizer $G_{S_l}$ of $S_l$ in $G$ on  ${\cal A}_{S_l}$ and on the category $\mbox{Nil}({\cal A}_{S_l})$ extend  to  actions of $\Sigma _{S_l} \ltimes \Bbb{Z}_p^n$.  By Proposition   \ref{changegroups}, we obtain a change of rings homomorphism $ \gamma: K^{\Sigma _{S_l} \ltimes \Bbb{Z}_p^n}({\cal A}_{S_l})  \rightarrow  K^{G _{S_l}}({\cal A}_{S_l})$ of $s\frak{I}$-algebras, as well as a homomorphism $\gamma ^{\mbox{Nil}}: K^{\Sigma _{S_l} \ltimes \Bbb{Z}_p^n}(\mbox{Nil}({\cal A}_{S_l}) ) \rightarrow K^{G _{S_l}}(\mbox{Nil}({\cal A}_{S_l}) )$ of $K^{\Sigma _{S_l} \ltimes \Bbb{Z}_p^n}({\cal A}_{S_l})$-modules.   Also by Proposition \ref{changegroups}, 
$K^{\Sigma _{S_l} \ltimes \Bbb{Z}_p^n}(\mbox{Nil}({\cal A}_{S_l}) )$ and $K^{G_{S_l}}(\mbox{Nil}({\cal A}_{S_l}) )$ are equipped with based structures over the $s\frak{I}$-algebras $K^{\Sigma _{S_l} \ltimes \Bbb{Z}_p^n}({\cal A}_{S_l}) $ and $K^{G _{S_l}}({\cal A}_{S_l}) $, respectively, and $\gamma ^{\mbox{Nil}}$ respects the based structures in the sense  of Proposition \ref{inflate}.  

\begin{Proposition}  With the data and assumptions of the preceding paragraph, the natural inflation map 
$$ K^{G _{S_l}}({\cal A}_{S_l}) \column{\wedge}{K^{\Sigma _{S_l} \ltimes \Bbb{Z}_p^n}({\cal A}_{S_l}) } K^{\Sigma _{S_l} \ltimes \Bbb{Z}_p^n}(\mbox{\em  Nil}({\cal A}_{S_l}) ) \rightarrow K^{G_{S_l}}(\mbox{\em  Nil}({\cal A}_{S_l}) )
$$
is an equivalence of $ K^{G _{S_l}}({\cal A}_{S_l})$-modules.  Consequently,  the set ${\cal S}(K^{G_{S_l}}(\mbox{\em  Nil}({\cal A}_{S_l}) ) = \{ \alpha _1, \alpha _2, \ldots \}$ has the property that $\alpha _i$ is in the image of the homomorphism  $\pi _0 (K^{\Sigma _{S_l} \ltimes \Bbb{Z}_p^n}({\cal A}^i_{S_l}) ) \rightarrow \pi _0(K^{G_{S_l}}({\cal A}^i_{S_l}) )$.  
\end{Proposition} 
\begin{Proof} The first follows directly from the discussion preceding the Proposition.  The second follows from statement (1) of Proposition \ref{msystemtech}. 
\end{Proof} 

\begin{Proposition}  The natural homomorphism 
$$R[\Sigma _{S_l} \ltimes \Bbb{Z}_p^n] \cong \pi _0 K^{\Sigma _{S_l} \ltimes \Bbb{Z}_p^n}(*) \rightarrow \pi _0 K^{\Sigma _{S_l} \ltimes \Bbb{Z}_p^n}({\cal A}^i_{S_l}) 
$$
is surjective.  
\end{Proposition}
\begin{Proof}  The group $\Sigma _{S_l} \ltimes \Bbb{Z}_p^n$ splits as  a product 
\begin{equation} \label{decomp} 
\Sigma _{S_l} \ltimes \Bbb{Z}_p^n \cong \Sigma _s \ltimes \Bbb{Z}_p^s \times \Sigma _t \ltimes \Bbb{Z}_p^t
\end{equation}
 where $s+t = n$.  Further, the special pro-scheme ${\cal A}$ splits as a product of special pro-schemes ${\cal X} $ and ${\cal Y}$,  where ${\cal X}$ is a product of $s$ copies of the system 
$$ \cdots \rightarrow   \Bbb{G}_m \stackrel{f^p}{\rightarrow}   \Bbb{G}_m \stackrel{f^p}{\rightarrow} \Bbb{G}_m \stackrel{f^p}{\rightarrow} \Bbb{G}_m $$
and ${\cal Y}$ is a product of $t$ copies of the system 
$$ \cdots \rightarrow \Bbb{A}_n  \stackrel{f^p}{\rightarrow}\Bbb{A}_n  \stackrel{f^p}{\rightarrow}\Bbb{A}_n \stackrel{f^p}{\rightarrow}\Bbb{A}_n  
$$
where in both cases the symmetric groups are acting by permuting coordinates and the groups $\Bbb{Z}_p^s$ and $\Bbb{Z}_p^t$ are acting by translations in the obvious way.  The origin in ${\cal Y}$ is an invariant subspace, and we have therefore have the inclusion ${\cal X} \times * \hookrightarrow {\cal X} \times {\cal Y}$, and it clearly induces an equivalence on equivariant $K$-theory by the homotopy property Proposition \ref{eqhomotopy}.  In terms of affine coordinate rings,  the action on ${\cal X} \times *$ can be identified with the action of $\Sigma _s \ltimes \Bbb{Z}_p^s \times \Sigma _t \ltimes \Bbb{Z}_p^t$ on the ring  $B_s = \bigcup _j k[x_1^{\pm \frac{1}{p^j}}, \ldots , x_s^{\pm \frac{1}{p^j}}]$, where the factor $\Sigma _s \ltimes \Bbb{Z}_p^s$ acts by permutations and translations in the obvious way, and where the factor $\Sigma _t \ltimes \Bbb{Z}_p^t$ acts trivially.  Galois theory of rings, as in \cite{magid}, now shows that the equivariant $K$-theory spectrum is given by 
$$  K^{\Sigma _{S_l} \ltimes \Bbb{Z}_p^n}({\cal X} \times {\cal Y}) \cong K(k[x_1^{\pm 1}, \ldots , x_s^{\pm 1}] \otimes _k k[\Sigma _t \ltimes \Bbb{Z}_p^t])
$$
The formula for the $K$-theory of Laurent extensions now shows that $\pi _0 K^{\Sigma _{S_l} \ltimes \Bbb{Z}_p^n}({\cal X} \times {\cal Y}) \cong R[\Sigma _t \ltimes \Bbb{Z}_p^t]$.  Because of the product decomposition (\ref{decomp}) above,  it is clear that the restriction homomorphism $R[\Sigma _{S_l} \ltimes \Bbb{Z}_p^n] \rightarrow R[\Sigma _t \ltimes \Bbb{Z}_p^t]$  is surjective, which is the required result.  
\end{Proof}  

\begin{Proposition}  Let $\{ \beta _i \} _{i \geq 1}$ denote the set ${\cal S}(K^{\Sigma _{S_l} \ltimes \Bbb{Z}_p^n}(\mbox{\em Nil} ({\cal A}_{S_l}))$, and let 
$$\nu _i : R[\Sigma _{S_l} \ltimes \Bbb{Z}_p^n] \cong \pi _0 K^{\Sigma _{S_l} \ltimes \Bbb{Z}_p^n}(*) \rightarrow \pi _0 
K^{\Sigma _{S_l} \ltimes \Bbb{Z}_p^n}( {\cal A}^i_{S_l})$$
denote that evident map induced by the map ${\cal A}^i_{S_l} \rightarrow *$. 
Let $B_i$ denote the $k$-algebra $k[x_1^{\frac{1}{p^{i}}}, \ldots , x_t^{\frac{1}{p^{i}}}]$ Let $\rho _i \in I_p[\Sigma _{S_l} \ltimes \Bbb{Z}_p^n] \subseteq R[\Sigma _{S_l} \ltimes \Bbb{Z}_p^n]$ denote the representation 
$$ k \column{\otimes}{B_{i-1}} B_i
$$ Then $\beta _i = \nu _i (\rho _i)$ for all $i$.  
\end{Proposition} 
\begin{Proof} The description of the elements $\beta _i$ is clear from the definition of the based structures on $K^G(\mbox{Nil}({\cal A}_{\frak{o}}))$ in Proposition \ref{indsystem}.  The fact that $\rho _i$ is in $I_p$ follows since the dimension of $\rho _i$ is $p^t$.  
\end{Proof}

\begin{Corollary}
Let $j^*:R[\Sigma _{S_l} \ltimes \Bbb{Z}_p^n] \rightarrow R[G_{S_l}]$ be the homomorphism induced by the representation $G_{S_l} \rightarrow \Sigma _{S_l} \ltimes \Bbb{Z}_p^n$.  Then the map 
$$ K^{G_{S_l}}(\mbox{\em Nil} ({\cal A}^i_{S_l})) \rightarrow  K^{G_{S_l}}(\mbox{\em Nil} ({\cal A}^{i+M}_{S_l}))
$$
in the direct system $K_{G_{S_l}}(\mbox{\em Nil} ({\cal A}_{S_l}))$ induces  multiplication by the element 
$j^*(\rho _i) j^*(\rho _{i+1}) \cdots j^*(\rho _{i+M}) $ in $R[G_{S_l}]$.  This element lies in $I_p[G_{S_l}]^M$.  
\end{Corollary} 
\begin{Proof} Clear. 
\end{Proof} 

\begin{Corollary}  There exists an integer so that the  map 
$$ K^{G_{S_l}}(\mbox{\em Nil} ({\cal A}^i_{S_l})) \rightarrow  K^{G_{S_l}}(\mbox{\em Nil} ({\cal A}^{i+M}_{S_l}))
$$
induces multiplication by an element of $I_p[G] \cdot R[G_{S_l}]$ for all $i$.  
\end{Corollary} 
\begin{Proof}  We note that we have a diagram
$$
\begin{diagram}
\node{G_{S_l}} \arrow{s} \arrow{e} \node{\Sigma _{S_l} \ltimes \Bbb{Z}_p^n} \arrow{s} \\
\node{G} \arrow{e} \node{\Sigma _n \ltimes \Bbb{Z}_p^n} 
\end{diagram} 
$$
and therefore a corresponding diagram of representation rings

\begin{equation}\label{repdiag}
\begin{diagram}
\node{R[\Sigma _n \ltimes \Bbb{Z}_p^n] } \arrow{s} \arrow{e} \node{R[\Sigma _{S_l} \ltimes \Bbb{Z}_p^n]} \arrow{s,t}{j^*} \\
\node{R[G]} \arrow{e} \node{R[G_{S_l}]} 
\end{diagram} 
\end{equation}
It follows from Corollary \ref{reptech} that there is an integer $M$ so that  $\rho _i \rho _{i+1} \cdots \rho _{i+M} \in I_p[\Sigma _n \ltimes \Bbb{Z}_p^n] \cdot R[\Sigma _{S_l} \ltimes \Bbb{Z}_p^n]$ for all $i$.  It follows that the map 
$$K^{G_{S_l}}(\mbox{Nil} ({\cal A}^i_{S_l})) \rightarrow  K^{G_{S_l}}(\mbox{Nil} ({\cal A}^{i+M}_{S_l}))
$$
induces multiplication by an element in $I_p[\Sigma _n \ltimes \Bbb{Z}_p^n] \cdot R[G_{S_l}]$ on homotopy groups  for all $i$.  Consequently, it induces multiplication by an element of $I_p[G] \cdot R[G_{S_l}]$ because of diagram (\ref{repdiag})  above.  
\end{Proof} 
\begin{Corollary}  We have 
$$ \Bbb{H}_p \column{\wedge}{ K^G(*)} (\column{hocolim}{\rightarrow} K^G(\mbox{\em Nil}({\cal A}_{\frak{o}}) )) \simeq *
$$
for all orbits $\frak{o}$ of the action of $G$ on the set ${\cal P}_i$ of subsets of cardinality $i$ of $\{ 1, \ldots , n \}$, where $ 1 \leq i \leq n-1$.  
\end{Corollary} 

\begin{Corollary} \label{affequiv}The natural map
$$   \Bbb{H}_p \column{\wedge}{K^G(*)}(\column{hocolim}{\rightarrow} K^G({\cal A};W))\rightarrow    \Bbb{H}_p \column{\wedge}{K^G(*)}(\column{hocolim}{\rightarrow }K^G({\cal G};W))
$$
is an equivalence of spectra.  It follows that the map
$$  \Bbb{H}_p \column{\wedge}{K^G(*)}(\column{hocolim}{\rightarrow} K^G(W))\rightarrow    \Bbb{H}_p \column{\wedge}{K^G(*)}(\column{hocolim}{\rightarrow }K^G({\cal G};W))
$$
induced by the projection ${\cal G} \times W \rightarrow W $ is an equivalence as well.  
\end{Corollary} 
\begin{Proof} The first part follows immediately from the facts that $$K^G(\mbox{Mod}_0({\cal A});W) \simeq K^G({\cal A};W)$$ and $$K^G(\mbox{Mod}_0({\cal G});W) \simeq K^G({\cal G};W)$$
The second statement is an immediate consequence of the homotopy property (Theorem \ref{eqhomotopy}).  
\end{Proof} 

\begin{Corollary}\label{products}  Let $\rho$ and $\rho^{\prime}$ denote monomial $\frak{G}$-representations of a pro-$p$ group $G$, and let    the corresponding pro-$G$ -schemes ${\cal G}$ and ${\cal G}^{\prime}$ denote the corresponding pro-$G$-schemes.  Then the natural map 
$$  \Bbb{H}_p \column{\wedge}{K^G(*)}(\column{hocolim}{\rightarrow} K^G({\cal G}))\rightarrow    \Bbb{H}_p \column{\wedge}{K^G(*)}(\column{hocolim}{\rightarrow }K^G({\cal G} \times{\cal G}^{\prime}))
$$
is an equivalence of spectra, where the map is the one induced by the projection ${\cal G} \times {\cal G}^{\prime} \rightarrow {\cal G}$.   
\end{Corollary} 
\begin{Proof} Easy extension of Corollary \ref{affequiv} to pro-schemes.  
\end{Proof} 
\section{The main theorem}

We are now in position to prove our main theorem.  

\begin{Theorem} \label{mainone} Let $G$ be a totally torsion free profinite $p$-group, and let $W$ be any Noetherian scheme with continuous $G$-action.  Then the natural map of $K^G(*)$-module spectra $K^G(W) \rightarrow K^G({\cal E} G \times W) $ induces an equivalence on completions 
$$ K^G(W)^{\wedge}_{\epsilon}  \rightarrow K^G({\cal E} G \times W) ^{\wedge}_{\epsilon}
$$
\end{Theorem} 
\begin{Proof} Since ${\cal E}G$ can be viewed as a pro object in the category of pro-schemes, in which every map is a projection ${\cal G} \times {\cal G}^{\prime} \rightarrow {\cal G}$, Corollary \ref{affequiv} tells us that the map 
$$  \Bbb{H}_p \column{\wedge}{K^G(*)} K^G(W) \rightarrow  \Bbb{H}_p \column{\wedge}{K^G(*)} K^G({\cal E}G;W)
$$  is an equivalence, it follows from Proposition \ref{criterion} that the map on completions is an equivalence.  \end{Proof} 

\begin{Corollary} \label{maintwo} Let $G$ and $W$  be as above, then there are equivalences of spectra 
$$  K^G(W)^{\wedge}_{\epsilon} \stackrel{\sim}{\rightarrow} K^G({\cal E}G \times W)^{\wedge}_p \cong K({\cal E}G \column{\times}{G} W) ^{\wedge}_p
$$
well defined up to homotopy. 
\end{Corollary}
\begin{Proof} Follows directly from Corollary \ref{nonequiv} and Proposition \ref{orbit}.  
\end{Proof}

 \end{document}